\documentclass[11pt,reqno,twoside,a4paper]{amsart}
\usepackage[latin1]{inputenc}
\usepackage[english]{babel}
\usepackage{graphicx,color,epsfig,latexsym,longtable,amsmath,amscd,amsxtra,amssymb,amsxtra,exscale,xy,mathrsfs,wrapfig,bbm}\xyoption{all}
\usepackage{multicol}
\theoremstyle{plain}
\newtheorem{theorem}{Theorem}[section]         
     
\newtheorem{lemma}[theorem]{Lemma}             
\newtheorem{proposition}[theorem]{Proposition} 
\theoremstyle{definition}
\newtheorem{example}[theorem]{Example}
\newtheorem{definition}[theorem]{Definition}   
               
\newtheorem{remark}[theorem]{Remark}

\renewcommand{\labelenumi} {(\alph{enumi})}    

\renewcommand{\theenumi} {(\alph{enumi})}

\begin{document}
  
\renewcommand {\Re} {\textit{Re}}
\renewcommand {\Im} {\textit{Im}}
\newcommand{\norm}[1] {\| #1 \|}
\newcommand{\lrnorm}[1]{\left\| #1 \right\|}
\newcommand{\bignorm}[1]{\bigl\| #1 \bigr\|}
\newcommand{\Bignorm}[1]{\Bigl\| #1 \Bigr\|}
\newcommand{\Biggnorm}[1]{\Biggl\| #1 \Biggr\|}
\newcommand{\biggnorm}[1]{\biggl\| #1 \biggr\|}
\newcommand {\sfrac}[2] { {{}^{#1}\!\!/\!{}_{#2}}} 
\newcommand {\einhalb} {\sfrac{1}{2}}
\newcommand {\pihalbe} {\sfrac{\pi\,}{2}}
\newcommand {\RR} {\mathbb R}
\newcommand {\NN} {\mathbb N}
\newcommand {\XX} {\mathbb X}
\newcommand {\ZZ} {\mathbb Z}
\newcommand {\CC} {\mathbb C}
\newcommand {\EE} {\mathbb E}
\newcommand {\FF} {\mathbb F}
\newcommand {\cA} {\mathcal A}
\newcommand {\cD} {\mathcal D}
\newcommand {\cF} {\mathcal F}
\newcommand {\cM} {\mathcal M}
\newcommand {\cP} {\mathcal P}
\newcommand {\cR} {\mathcal R}
\newcommand {\cT} {\mathcal T}
\newcommand {\cS} {\mathcal S}
\newcommand {\cU} {\mathcal U}
\newcommand {\cX} {\mathcal X}
\newcommand {\cY} {\mathcal Y}
  \newcommand {\eins} {\mathbbm 1}
  \newcommand {\la}{\lambda}
  \newcommand {\si}{\sigma}
  \newcommand {\Si}{\Sigma}
  \newcommand {\Ga}{\Gamma}
  \newcommand {\ga}{\gamma}
  \newcommand {\om}{\omega}
  \newcommand {\Om}{\Omega}
  \newcommand {\KN}[1] {|\!|\!|#1|\!|\!|}
  \newcommand {\al}{\alpha}
  \newcommand {\eps}{\epsilon}
  \newcommand {\Control}[1]{\Phi_{#1}}
  \newcommand {\Obs}[1] {\Psi_{#1}}
  \newcommand {\IO}[1] {\FF_{#1}}
  \newcommand {\ControlMap}{\Control{\tau}}
  \newcommand {\ObsMap}{\Obs{\tau}}
  \newcommand {\ECM}{\Control{\infty}}
  \newcommand {\EOM}{\Obs{\infty}}
  \newcommand {\IOM}{\IO{\infty}}
  \newcommand {\Sec}[1] {S({#1})}
  \newcommand {\Continuous} {\sC}
  \newcommand {\ran}[1] {\textit{ran}\,(#1)}
  \newcommand {\bigidual}[3] {\bigl\langle #1, #2 \bigr\rangle_{#3}}
  \newcommand {\idual}[3] {\langle #1, #2 \rangle_{#3}}
  \newcommand {\bigdual}[2] {\bigidual{#1}{#2}{}}
  \newcommand {\dual}[2] {\idual{#1}{#2}{} }
  \newcommand {\ws} {weak$^*$}
  \newcommand {\weg} {\backslash}
  \newcommand {\bbrace}[1] {[\![ #1 ]\!]} 
  \newcommand {\emb} {\hookrightarrow}
  \newcommand {\Rad} {\text{Rad}}
\newcommand {\rA} {(\la{+}A)^{-1}}
\newcommand {\rAm} {(\la{+}\Am)^{-1}}
\newcommand {\Xpe} {\dot{X}_1}
\newcommand {\Xpm} {\dot{X}_{-1}} 
\newcommand {\Xpepm} {(\Xpe)^{\dot{}}_{-1}}
\newcommand {\Xpmpe} {(\Xpm)^{\dot{}}_{1}}
\newcommand {\Ape} {\dot{A}_1}
\newcommand {\Apm} {\dot{A}_{-1}} 
\newcommand {\Xm}{{X_{-1}}}
\newcommand {\Tm}{{T_{-1}}}
\newcommand {\Am}{{A_{-1}}}
\newcommand {\Amm}{{(A_{-1})^{-1}}}
\newcommand {\McIntosh} {{M\textsuperscript{c}Intosh}{ }}
\newcommand {\nnorm}[1] {|\!|\!|#1|\!|\!|}
\newcommand {\LpAbsch}[1]{${L_\ast^{#1}}$}
\def\eqnaintertext#1{\crcr\noalign{\vskip\belowdisplayskip
 \vbox{\normalbaselines\noindent#1}\vskip\abovedisplayskip}}
\newcommand{\Id}{\text{Id}}
\renewcommand{\thefootnote}{\fnsymbol{footnote}}

  \date{18th July 2005 (Revision 30.3.2006)}
  \title[Weighted Admissibility and Wellposedness in Banach spaces] 
        {Weighted Admissibility and Wellposedness of linear
          systems in Banach spaces} 
  
  \author[Bernhard H. Haak, Peer Chr. Kunstmann]{Bernhard H. Haak
    ${}^1$ and Peer Chr. Kunstmann ${}^1$}
  \address{Mathematisches Institut I\\ Universität Karlsruhe\\
           Englerstraße 2\\ 76128 Karlsruhe\\ Germany}       
  \email{Bernhard.Haak@math.uni-karlsruhe.de \\
         Peer.Kunstmann@math.uni-karlsruhe.de}
  \thanks{${}^1$The research is supported in part by the DFG project
  \glqq{}$H^\infty$--Kalk\"ul und seine Anwendungen auf partielle
  Differentialgleichungen\grqq{} under contract number WE 2847/1-1.}
\begin{abstract}
We study linear control systems in infinite--dimensional Banach spaces
governed by analytic semigroups.
For $p\in[1,\infty]$ and $\al\in\RR$ we introduce the notion of 
$L^p$--admissibility of type $\al$ for unbounded observation and 
control operators. Generalising earlier work by Le~Merdy 
\cite{LeMerdy:weiss-conj} and the first named 
author and Le~Merdy \cite{HaakLeMerdy} we give conditions under 
which $L^p$--admissibility of type $\al$ is characterised by 
boundedness conditions which are similar to those in the well--known 
Weiss conjecture. We also study $L^p$--wellposedness of type $\al$
for the full system. Here we use recent ideas due to Pruess and 
Simonett. Our results are illustrated by a controlled heat equation
with boundary control and boundary observation where we take Lebesgue
and Besov spaces as state space. This extends the considerations in
\cite{BGSW} to non--Hilbertian settings and to $p\neq 2$. 
\end{abstract}
  \maketitle

\section{Introduction and Main Theorems}\label{sec:intro}
We are concerned with linear control systems of the following form
\begin{equation}\label{eq:control-system}
\left\{
\begin{array}{lcl}
 x'(t)+A x(t) & = & B u(t),\quad (t>0)\\
 x(0) & = & x_0,\\
 y(t) & = & C x(t),\quad (t>0),
\end{array}
\right.
\end{equation}
where $-A$ is the generator of a strongly continuous semigroup
$(T(\cdot))$ in a Banach space $X$. The function $x(\cdot)$ takes
values in $X$, the functions $u(\cdot)$ and $y(\cdot)$ take values in
Banach spaces $U$ and $Y$, respectively.  The control operator $B$ is
an unbounded operator from $U$ to $X$, and the observation operator
$C$ is an unbounded operator form $X$ to $Y$. We refer e.g. to
\cite{Salamon, Staffans,Weiss:admiss-observation, Weiss:Admissibility-of-unbounded}.
A commonly used minimal assumption on $B$ and $C$ is
that $C$ is  bounded $X_1(A)\to Y$ and $B$ is bounded $U\to X_{-1}(A)$
where $X_1(A)$ denotes  the domain $D(A)$ of $A$ equipped with the
graph norm, and $X_{-1}(A)$ denotes  the completion of $(X,\norm{ 
R(\la_0,A)\cdot}_X)$ with $\la_0$ in the resolvent set $\rho(A)$ of
$A$ (cf., e.g., \cite[Sect. II.5]{EngelNagel}).  Note that, for
$\la_0\in\rho(A)$, the norm $\|(\la_0-A)\cdot\|_X$ on $X_1(A)$ is
equivalent to the graph norm of $A$. The semigroup $T(\cdot)$ has an
extension to a strongly continuous semigroup $\Tm(\cdot)$ on
$\Xm:=X_{-1}(A)$ whose generator $\Am$ is an extension of $A$
(cf. \cite[Sect. II.5]{EngelNagel}). These extensions are needed in order to
give a precise meaning to compositions involving the operator $B$.
Now let $\cX :=C([0,\infty),X)$ and, for each $\tau>0$,
$\cX_\tau=C([0,\tau],X)$. Suppose that we are given spaces 
$\cY$ of functions $\RR_+\to Y$ and $\cU$ of functions
$\RR_+\to U$ with restrictions $\cY_\tau$, $\cU_\tau$
to $[0,\tau]$, $\tau>0$, respectively. Then the system
\eqref{eq:control-system} is called \emph{wellposed}  if, for each
$\tau>0$, state and output of the system \eqref{eq:control-system}
depend continuously on initial state and input, i.e., if the mapping
\[
  X\times\cU_\tau \to X\times\cY_\tau, \quad
  (x_0,u(\cdot))\mapsto (x(\tau),y(\cdot))
\]
is continuous. Since the solution to \eqref{eq:control-system} has the
(formal) representation
\[
 \bigl(x(t),y(t)\bigr)=\bigl( T(t)x_0+\int_0^t \!\!T(s)Bu(t-s)\,ds, \;
              CT(t)x_0+C\!\!\int_0^t \!\!T(s)Bu(t-s)\,ds\bigr),\quad t>0,
\]
and $T(\cdot)$ is strongly continuous, this means continuity of the three 
maps 
\begin{equation*} \label{def:Phi-Psi-cF}
\Psi_\tau: \left\{\begin{array}{l}
 X \to    \cY_\tau \\
x_0\mapsto CT(\cdot)x_0
\end{array}\right.\!\!, \quad
\Phi_\tau: \left\{\begin{array}{l}
\cU_\tau \to     X\\
u                \mapsto \int\limits_0^\tau \Tm(\tau{-}s)B u(s)\,ds 
\end{array}\right. \text{ and }\;
\cF_\tau:
\left\{\begin{array}{l}
 \cU_\tau\to      \cY_\tau,\\
 u               \mapsto  C \Tm(\cdot)B * u
\end{array}\right.
\end{equation*}
where $*$ denotes convolution. Requiring continuity of these maps
leads to the notions of \emph{admissibility} for the observation
operator $C$, for the control operator $B$, respectively, and the
notion of {\em wellposedness} for the input--output map
$\cF_\tau$ of system \eqref{eq:control-system}.  
In case of uniformly exponentially stable semigroups and
$\cY_\tau=L^p([0,\tau],Y)$ and $\cU_\tau=L^p([0,\tau],U)$, one may also
consider wellposedness on $\RR_+=(0,\infty)$, i.e. the spaces 
$\cY=L^p(\RR_+,Y)$ and $\cU=L^p(\RR_+,U)$ and continuity of the maps
\[
\Psi: \left\{\begin{array}{l}
 X \to     \cY \\
x_0\mapsto CT(\cdot)x_0
\end{array}\right.,\qquad
\Phi: \left\{\begin{array}{l}
\cU_\tau \to     X\\
u                \mapsto \int\limits_0^\infty \Tm(t)B u(t)\,dt 
\end{array}\right. \text{ and }\;
\cF:
\left\{\begin{array}{l}
 \cU \to     \cY,\\
 u           \mapsto  C\Tm(\cdot)B *u
\end{array}\right..
\] 
Notice that in the second mapping the convolution is replaced by an
integration of the semigroup against the right hand side
(cf. \cite{Weiss:Admissibility-of-unbounded} and
Remark~\ref{rem:dualisieren-von-lp-al-zulaessig} below). 
\smallskip
If $X$, $Y$ and $U$ are Hilbert spaces it is natural to take 
$\cY=L^2(\RR_+,Y)$ and $\cU=L^2(\RR_+,U)$
(cf.\cite{JacobPartington:survey,Weiss:Admissibility-of-unbounded,Weiss:admiss-observation}). 
In the Banach space case one can consider $\cY=L^p(\RR_+,Y)$ and 
$\cU=L^p(\RR_+,U)$ where $p\in[1,\infty]$
(cf. \cite{Engel:GrabowskiCallier,JacobPartington:survey,Staffans:Lp-multipliers,Weiss:Admissibility-of-unbounded}).
For the case $p=2$, the notion of \emph{$\al$--admissibility} for observation 
and control operators was introduced in \cite{HaakLeMerdy} meaning that the 
$L^2$--space on $\RR_+$ with values in $Y$ or $U$, respectively, is taken with 
respect to a polynomial weight on $\RR_+$. In this paper, we will extend that 
notion for observation operators from $L^2$--norms to $L^p$--norms
with $p\in[1,\infty]$. For control operators our definition differs
from that given in \cite{HaakLeMerdy}
(cf. Remarks \ref{rem:zusammenhang-HaLeM} and
\ref{rem:dualisieren-von-lp-al-zulaessig} below). Furthermore, we
will introduce and study the new notion of \emph{$L^p$--wellposedness of type 
$\alpha$} for the system \eqref{eq:control-system}.
For a short discussion on our motivation we refer to Theorem~\ref{thm:nonlin}
and Remark~\ref{rem:motivation}.
The main basic question in modelling a given system in order to obtain a wellposed 
system of the form \eqref{eq:control-system} is, of course, how to check for
admissibility of the operators $C$ and $B$ and the wellposedness of the
input--output map $\cF$.
For Hilbert spaces $X$, $Y$, $U$, the well-known \emph{Weiss conjecture}
(cf. \cite{JacobPartington:survey,Weiss:conjectures}) relates
$L^2$--admissibility on $\RR_+$ of $C$ and $B$ to boundedness of
\begin{eqnarray*}
 W_C&:=&\{ \la^\einhalb C(\la+A)^{-1}:\; \la>0\}\subseteq B(X,Y)
 \qquad\mbox{and}\\
 W_B&:=&\{ \la^\einhalb (\la+\Am)^{-1}B:\; \la>0\}\subseteq B(U,X),
\end{eqnarray*}
respectively. Assuming that $T(\cdot)$ is bounded analytic and that
$A^\einhalb$ is  $L^2$--admissible for $A$, boundedness of $W_C$
actually characterises  $L^2$--admissibility of $C$
(cf. \cite{LeMerdy:weiss-conj} for this result and a  detailed
discussion of the Weiss conjecture).  In \cite{HaakLeMerdy}, this
characterisation was extended to cover the case of
$L^2$--admissibility of type $\al$ for a certain range of $\al$.
\smallskip 
For $\om \in(0,\pi)$, let  $S(\om) := \{z\in\CC \setminus\{0\}:|\arg z| < \om\}$
and let $S(0) := (0,\infty)$. We recall that a \emph{sectorial
  operator $A$ of type $\om\in[0,\pi)$} in a Banach space $X$ is a
closed linear operator $A$ satisfying $\sigma(A)\subseteq
\overline{S(\om)}$ and, for any $\nu\in(\om,\pi)$, 
\[
 \sup\{ \|\la R(\la,A)\|: |\arg\la|\ge\nu\}<\infty.
\]
Observe that $-A$ generates a bounded analytic semigroup in $X$ if and only if 
$A$ is a densely defined sectorial operator in $X$ of type $<\pihalbe$. 
Moreover, we recall that a sectorial operator with dense range is actually 
injective (cf. \cite[Thm. 3.8]{McIntoshYagi:without}). Occasionally we
shall use fractional powers of sectorial operators or their functional
calculus. We refer to, e.g., \cite{Haase:Buch, KunstmannWeis:Levico, McIntosh:H-infty-calc}.
\subsection*{$L^p$--admissibility of type $\alpha$}\label{sec:Lp-zul-typ-alpha}
In this paper, we assume that $A$ is a densely defined sectorial operator 
of type $<\pihalbe$ with dense range, and we characterise $L^p$--admissibility 
of type $\al$ for observation  and control operators, the range of values 
of $\al$ depending on $p\in[1,\infty]$. For a fixed
$p\in[1,\infty]$ we denote by $p'$ the dual exponent given by
$\sfrac{1}{p}+\sfrac{1}{p'}=1$. In order to formulate our first
definition we  introduce, for any $k\in\NN$ and a given generator $-A$
of a bounded strongly continuous semigroup $T(\cdot)$ in $X$, the spaces
$X_k:=X_k(A)$, which is $\cD(A^k)$  equipped with the norm
$\|(\text{Id}+A)^k\cdot\|_X$, and $X_{-k}:=X_{-k}(A)$, which is  the
completion of $X$ with respect to the norm $\|(\text{Id}+A)^{-k}\cdot\|_X$.  
Observe that $1\in\rho(-A)$ and that replacing $\text{Id}+A$ by $\la_0+A$ in
these  expressions leads to equivalent norms whenever
$\la_0\in\rho(-A)$  (again, we refer to \cite[Sect. II.5]{EngelNagel}).
The semigroup $T(\cdot)$ has an extension to a bounded strongly continuous
semigroup $T_{-k}(\cdot)$ on $X_{-k}$ whose generator $A_{-k}$ is an
extension of the  operator $A$.
\medskip
We denote, for $\al\in\RR$ and intervals $I\subseteq \RR_+$,  
\[
L_\al^p(I, X) := \bigl\{ f: I \to X:\; t\mapsto t^\al f(t) \in
L^p(I,X) \bigr\}.
\]
The space $L^p_\al(\RR_+, X)$ is abbreviated $L^p_\al(X)$.

\begin{definition}\label{def:Lp-al-zul}
Let $X$, $U$, $Y$ be Banach spaces and for some $k\in\NN$ let 
$C \in B(X_k,Y)$ and $B\in B(U, X_{-k})$. 
Given $p\in[1,\infty]$ and a bounded analytic semigroup $T(\cdot)$
on $X$,
\begin{enumerate}
\item \label{item:def-C-Lp-al-zul}
  $C$ is called {\em finite--time $L^p$--admissible of type}
  $\al>-\sfrac{1}p$, if for any $\tau>0$ there is a constant
  $M_\tau>0$ such that for all $x\in X_k$ we have
  \begin{equation}\label{eq:def-Lp-al-zul-obs-finite}
   \bignorm{ t \mapsto C T(t)x }_{L^p_\al([0,\tau],Y)} \le M_\tau \norm{x}_X.
  \end{equation}
\item \label{item:def-B-Lp-al-zul}
  $B$ is called {\em finite--time $L^p$--admissible of type}
  $\al < \sfrac{1}{p'}$ (or $\al\le 0$ for $p=1$), if for any $\tau>0$
  there is a constant $K_\tau>0$ such that for all  $u\in
  L^p_\al((0,\tau),U)$ we have
  \begin{equation}\label{eq:def-Lp-al-zul-contr-finite}
      \bignorm{T_{-k}B * u }_{L^\infty((0,\tau), X )}
        \le K_\tau \norm{u}_{L^p_\al((0,\tau),U)}.
  \end{equation}
\end{enumerate}
If the above estimates hold with constants $M_\tau$ and $K_\tau$ that
can be chosen independently of $\tau>0$, the operators $B$ and $C$ are
called (infinite--time) {\em $L^p$--admissible of type $\al$}.
\end{definition}

\begin{remark}\label{rem:zusammenhang-HaLeM}
In part \ref{item:def-B-Lp-al-zul} of the definition, 
inequality \eqref{eq:def-Lp-al-zul-contr-finite} means
ess$\,\sup_{t\in(0,\tau)} \norm{ \Phi_t( u ) }_X \le
K_\tau \norm{ u }_{ L^p_\al((0,\tau),U)}$. This seems to be better
suited when dealing with weighted spaces than requiring just 
$\norm{ \Phi_\tau(u) }_X \le K_\tau\norm{ u }_{L^p_\al((0,\tau),U)}$.
\smallskip
Notice that our definition differs from that in \cite{HaakLeMerdy}:
for observation operators, $L^2$--admissibility of type $\al$ is
called $2\al$--admissibility there. For control operators, our
definition is different from the notion studied in
\cite{HaakLeMerdy}. This is due to the application to nonlinear
problems in Theorem~\ref{thm:nonlin}. See also
Remark~\ref{rem:dualisieren-von-lp-al-zulaessig} for
a study of both notions. 
\end{remark}
\smallskip

\begin{lemma}\label{lem:alternative-bed-fuer-obs-zul}
The notion of finite--time $L^p$--admissibility of type $\al$ for
$A$ is independent of the underlying interval $[0,\tau]$ in the definition. 
If the semigroup is uniformly exponentially stable and $\al \ge 0$ or 
$p=\infty$, finite--time $L^p$--admissibility of type $\al$ for control
operators is equivalent to (infinite--time) $L^p$--admissibility of
type $\al$. For observation operators this is true for all
$\al$. Moreover, an observation operator $C$ is (infinite--time)
$L^p$--admissible of type $\al$ of $A$ if and only if 
\begin{equation} \label{eq:def-Lp-al-zul-obs}
  \biggl( \int_0^\infty \norm{ t^\al CT(t) x }^p\,dt \biggr)^{\sfrac{1}p}
  \le M_{\tau,\al} \norm{ x }.
\end{equation}
holds. 
\end{lemma}
Notice that the requirement $\al \ge 0$ or $p=\infty$ in the above
equivalence assertion for control operators is necessary. We provide a
short counterexample in case $\al < 0$ in
Example~\ref{ex:gegenbeispiel-alpha-negativ}.
\begin{remark}[Dualisation]\label{rem:dualisieren-von-lp-al-zulaessig}
Let $X$, $Y$ be reflexive Banach spaces and $p\in (1,\infty)$ and
assume that $C \in B(X_1,Y)$ is $L^p$--admissible of type $\al$ for
$A$ on $\RR_+$, Then, for $u \in L_{-\al}^{p'}(\RR_+, Y')$, 
\begin{eqnarray*}
  \bigidual{t\mapsto CT(t)x}{u} {L_\al^p \times L_{-\al}^{p'}}
&=& \int_0^\infty \bigidual{t^\al CT(t)x}{t^{-\al} u(t)} {Y \times Y'}\, dt\\
&=& \bigidual{x}{\int_0^\infty  T(t)' C'  u(t) \, dt} {X \times X'}.
\end{eqnarray*}
One can therefore consider the following dual condition to
(\ref{eq:def-Lp-al-zul-obs}), that was introduced in
\cite{HaakLeMerdy}.
\begin{equation} \label{eq:def-Lp-zul-contr-falsch}
  \biggnorm{ \int_0^\infty T(t) B u(t) \,dt }_X \le K
  \norm{u}_{L^p_{-\al}(\RR_+, U)},
\end{equation}
where the integral is considered as a Pettis integral in $X_{-k}$
taking values in $X$. Notice, that in case $\al\not=0$ the reflection
$R_\tau u = u(\tau{-}\cdot)$ is not bounded on $L_\al^p([0,\tau])$. 
Therefore, one cannot hope that (\ref{eq:def-Lp-zul-contr-falsch})
might be equivalent to infinite--time $L^p$--admissibility of type $\al$
unless $\al=0$. We shall see in Theorem~\ref{thm:control-Lp-zul-typ-alpha}
below that condition (\ref{eq:def-Lp-zul-contr-falsch}) is indeed a
stronger notion than $L^p$--admissibility of type $\al$.
However, for $\al>0$, $L^p$--admissibility of
an observation operator $C$ {\em implies} $L^p$--admissibility of $C'$
on $X'$ by Theorem~\ref{thm:control-Lp-zul-typ-alpha} below.
\end{remark}
\medskip

\subsection*{$L^p_\ast$--estimates and the real interpolation method}
As mentioned above, the crucial condition in \cite{LeMerdy:weiss-conj} 
(besides $T(\cdot)$ being bounded and analytic) for the Weiss
conjecture to hold was that $A^\einhalb$ is admissible for $A$, i.e.,
the existence of a constant $L>0$ such that 
\[
 \int_0^\infty \bignorm{ A^\einhalb T(t)x}^2 \,dt
= 
 \int_0^\infty \bignorm{(tA)^\einhalb T(t)x}^2 \,\tfrac{dt}t
=
 \bignorm{\psi(tA) x}_{L^2(\RR_+, dt/t, X)}^2 \le L\norm{ x }_X^2,
\]
where $\psi(z)=z^\einhalb e^{-z}$. In our situation this corresponds to
\begin{equation}
  \label{eq:Lp-stern-abschaetzungen-und-zul}
 \int_0^\infty \bignorm{t^\al A^{\al+\sfrac{1}{p}} T(t)x}^p \,dt
= 
 \int_0^\infty \bignorm{(tA)^{\al+\sfrac{1}{p}} T(t)x}^p \,\tfrac{dt}t
=
\bignorm{\psi(tA) x}_{L^p(\RR_+, dt/t, X)}^p\le \widetilde L \norm{ x }_X^p
\end{equation}
for $\psi(z) = z^{\al+\sfrac{1}{p}} \exp(-z)$. 
We denote $L^p_*(\RR_+, X) := L^p(\RR_+, dt/t, X)$.
It is known that the property of a sectorial operator $A$ of type $\om$ on a Banach
space $X$, to satisfy an estimate 
$\norm{ \psi(\cdot A)x }_{L^p_*(\RR_+,Y)}\le L\norm{ x }_X$ does not depend upon the 
particular choice of the function $\psi\in H_0^\infty(\Sec{\nu})\weg\{0\}$, 
$\nu>\om$, where
\[
H_0^\infty(\Sec{\nu}) = \bigl\{ f\in H^\infty(\Sec{\nu}):\; 
\exists c,s>0:\; \bigl|f(z)\bigr| \le c \,\tfrac{|z|^s}{1+|z|^{2s}} \bigr\},
\]
see \cite[Thm. 5]{McIntoshYagi:without},
\cite[Thm. 4.1]{AuscherMcIntoschNahmod:holomorphic},
\cite[Thm. 4.3]{Haase:H-infty-calc} or
\cite[Thm 6.4.3]{Haase:Buch}. Therefore, we say that 
{\em $A$ satisfies \LpAbsch{p}--estimates on $X$} in this case.
\medskip
Let $A$ be a densely defined sectorial operator of type $\om$ with
dense range on $X$. Let $\XX := (X, \norm{A(I+A)^{-2}  \cdot})^\sim$. 
For some holomorphic function $\psi$ on $\Sec{\nu}$,
$\nu>\om$ and $\theta\in (-1,1)$ such that $z^{-\theta}\psi(z) \in
H_0^\infty(\Sec{\nu})\weg\{0\}$ consider the space 
\[
X_{\theta,\psi,p} := \left\{ x\in \XX: \; t^{-\theta} \psi(tA)x \in
  L^p_\ast(\RR_+, X) \right\}.
\]
Since $X_{\theta,\psi,p}$ does not depend on $\psi \in
H_0^\infty(\Sec{\nu})\weg\{0\}$ (see above), we write for short $X_{\theta,p}
:= X_{\theta,\psi,p}$.
Resorting to the space $\XX$ in this setting allows explicitly that
$X_{\theta,p}$ may be a larger space than $X$. For the background of
this construction we refer e.g. to 
\cite{KaltonWeis:euclidian-structures}, \cite[Section 2]{HaakHaaseKunstmann} or
\cite[Chapter 6.3]{Haase:Buch}.
Notice that in this terminology \LpAbsch{p}--estimates for $A$ read as
$X \emb X_{0,p}$. These spaces are strongly connected to
real interpolation spaces between the homogeneous spaces
$\dot{X}_{-1}$ and $\dot{X}_1$.
Here, for $\sigma\in\ZZ$, we denote by $\dot{X}_\sigma:=\dot{X}_{\sigma}(A)$
the completion of the space $\cD(A^\sigma)$ with respect to the homogeneous
norm $\norm{ A^\sigma\cdot }_X$. We refer to \cite{KaltonKunstmannWeis},
\cite[Sect. 15]{KunstmannWeis:Levico}, \cite{HaakHaaseKunstmann} 
or \cite{Haase:H-infty-calc} for more details on these spaces. Observe
that -- due to our  assumptions -- the space $\cD(A^\sigma)$ is dense in
$X$, and that  $\dot{X}_\sigma=X_\sigma$ in case $0\in\rho(A)$. With
these notations the following result holds true
(cf. \cite{AuscherMcIntoschNahmod:holomorphic} for the case $p=2$
and \cite[Thm. 6.4.6]{Haase:Buch}, \cite[Thm. 5.2]{Haase:H-infty-calc}
for the general case).
\begin{theorem}\label{thm:char-McI-Y-spaces}
Let $X$ be a Banach space  and $A$ be a densely defined sectorial
operator of type $\om$ with dense range on $X$. Then
\[
\bigl(\dot{X}_{-1}(A), \dot{X}_1(A))_{\theta,p} = X_{2\theta-1, p}
\]
for all $\theta\in(0,1)$ and $p\in[1,\infty]$.
\end{theorem}
\begin{remark}\label{rem:Lp-Absch-chat-real-interpol}
In the situation of the above theorem, 
$A$ has $L^p_\ast$--estimates if and only if 
$X \emb \bigl(\dot{X}_{-1}(A), \dot{X}_1(A))_{\einhalb,p}$.
This can be used to establish $L^p_\ast$--estimates for given operators
in concrete cases, as we will show in Section \ref{sec:appl}.
\end{remark}
We start our main results with the following characterisations of
$L^p$--admissibility of type $\al$ for observation and control
operators that extend the results in \cite{HaakLeMerdy}.
\begin{theorem}\label{thm:obs-Lp-zul-typ-alpha}
Let $p\in[1,\infty]$ and $A$ be a densely defined sectorial operator
of type $\om<\pihalbe$ with dense range on $X$.
Let $C \in B(X_k, Y)$ be an observation operator for some $k\geq
1$. Let $\alpha\in(-\sfrac{1}p, k-\sfrac{1}{p})$ and consider
the set 
\begin{equation}\label{eq:W_C}
W_C := \bigl\{ t^{k-\alpha-\sfrac{1}{p}}\,C (t+A)^{-k}: \;
t>0\bigr\} \subseteq B(X,Y).
\end{equation}
Then the following assertions hold:
\begin{enumerate}
\item If $C$ is $L^p$--admissible of type $\alpha$ for $A$, then $W_C$
  is bounded in $B(X,Y)$.
\item If $A$ satisfies \LpAbsch{p}--estimates and $W_C$ is bounded in $B(X,Y)$,
  then $C$ is $L^p$--admissible of type $\alpha$ for $A$.
\end{enumerate}
\end{theorem}
Notice that, for the operator $C := A^{\al+{\sfrac{1}p}}$, the set
$W_C$ in \eqref{eq:W_C} is bounded, whence the assumption of
\LpAbsch{p}--estimates in the above characterisation
cannot be improved (cf. \eqref{eq:Lp-stern-abschaetzungen-und-zul}). 
\begin{theorem}\label{thm:control-Lp-zul-typ-alpha} 
Let $p\in[1,\infty]$ and $A$ be a densely defined sectorial operator
of type $\om<\pihalbe$ with dense range on $X$. 
Let $B \in B(U,X_{-k})$ be a control operator for some $k\geq 1$ and
let $\alpha\in(\sfrac{1}{p'}-k, \sfrac{1}{p'})$ and consider the set
\begin{equation}\label{eq:W_B}
W_B := \bigl\{ t^{k+\alpha-\sfrac{1}{p'}}\,(t+A_{-k})^{-k}B: \; t>0 \bigr\}
       \subseteq B(U,X).
\end{equation}
Then the following assertions hold:
\let\ALTA\labelenumi \let\ALTB\theenumi
\def\theenumi{(\alph{enumi})$\;$} \def\labelenumi{(\alph{enumi})}
\begin{enumerate}
\item \label{item:thm-control-zul-equiv-notw}
      If $B$ is $L^p$--admissible of type $\alpha$ for $A$, then $W_B$
      is bounded.
\item \label{item:thm-control-zul-equiv-hinr}
      If $W_B$ is bounded and $\al>0$ (for $p>1$) or $\al=0$ (for
      $p=1$) then $B$ is $L^p$--admissible of type $\alpha$ for $A$.
\end{enumerate}
\def\theenumi{(\alph{enumi})'} \let\labelenumi\theenumi
\begin{enumerate}
\item \label{item:thm-control-duale-bed-notw}
      If the dual condition (\ref{eq:def-Lp-zul-contr-falsch}) holds,
      then $W_B$ is bounded.
\item \label{item:thm-control-duale-bed-hinr}
      If $p<\infty$ and $W_B$ is bounded and the adjoint operator $A'$ satisfies
      \LpAbsch{p'}--estimates, then (\ref{eq:def-Lp-zul-contr-falsch})
      holds.
\end{enumerate}
\let\labelenumi\ALTA \let\theenumi\ALTB
\end{theorem}
As mentioned above, in case $\al=0$, $L^p$--admissibility of type
$\al$ and (\ref{eq:def-Lp-zul-contr-falsch}) are
equivalent. In particular, a modified Weiss conjecture on control
operators holds for $p=1$.
\bigskip
Also the boundedness conditions on the sets $W_C$ and $W_B$ in
(\ref{eq:W_C}) and (\ref{eq:W_B}) may be characterised by real interpolation
methods. In fact, they are conditions on the {\em domain} of (a
continuous extension of)  the observation operator $C$ and the {\em
  range} of the control  operator $B$, respectively. This way of
viewing the boundedness conditions on $W_C$ and $W_B$, respectively,
is not new for control operators (cf. \cite{Weiss:conjectures}) under
the supplementary assumption of $0\in\varrho(A)$. The general
equivalence for observation and control operators was first studied in
\cite{HaakHaaseKunstmann}. The proof of the following is very similar 
to the arguments used there and thus omitted.
\begin{theorem}\label{thm:beschr-bedingung-und-interpol}
Let $A$ be an injective  sectorial operator on the Banach space
$X$, and let $k\in\NN$. Let $C \in B(X_k,Y)$ and $B \in B(U,X_{-k})$ be
bounded operators,  where $U$, $Y$ are Banach spaces, and let 
$\theta \in (0,1)$. Then the following equivalences hold true:
\begin{enumerate}
\item \label{item:beschr-bedingung-und-interpol-a}
  $\sup\limits_{0 < \lambda < \infty} \norm{\lambda^{k\theta} 
  (\lambda+ A_{-k})^{-k} B}_{U\to X} < \infty$  if and only if
  $B: U \to X_{-k}$ is bounded in the norm $U \to (\dot{X}_{-k}, X)_{\theta,\infty}$.
\item \label{item:beschr-bedingung-und-interpol-b}
  $\sup\limits_{0 < \lambda < \infty} \norm{\lambda^{k(1-\theta)}C (\lambda+
  A)^{-k}}_{X\to Y} < \infty$  if and only if
  $C: X_k \to Y$ is bounded in the norm $(X, \dot{X}_k)_{\theta,1} \to Y$.
\end{enumerate}
\end{theorem}
\begin{remark}\label{rem:k-egal}
Of course, if for some $k\in\NN$, $C \in B(X_k, Y)$ and if the
operators $\lambda^{k(1-\theta)}C (\lambda+ A)^{-k}$, $\la>0$, are
uniformly bounded, then for any $n\in\NN_0$, $C \in B(X_{k+n}, Y)$
and by the sectoriality of $A$ the operators $\lambda^{n+k(1-\theta)}C (\lambda+ A)^{-k-n}$,
$\la>0$, are uniformly bounded. However, by
reiteration (cf. \cite[1.10.2]{Triebel:interpolation}), 
\[
\bigl(X, \dot{X}_k\bigr)_{\theta,1} = \bigl(X, \dot{X}_{k+n} \bigr)_{\si,1}
\]
for $\si = \tfrac{n+k(1-\theta)}{k+n}$, whence the assertions of 
Theorem~\ref{thm:beschr-bedingung-und-interpol} do not
depend on the question, for which $k\in\NN$ the given boundedness conditions are
satisfied. The same arguments apply to the conditions on control operators.
\end{remark}
\medskip
\subsection*{Wellposedness of the full system }
We now study the full system (\ref{eq:control-system}) for 
$\cY = L_\al^p(\RR_+, Y)$ and $\cU = L_{\al}^p(\RR_+, U)$. 
\begin{definition}
Let $X, U, Y$ be Banach spaces and let $k\in\NN$, $p\in[1,\infty]$
and $\al\in(-\sfrac{1}p, \sfrac{1}{p'})$. 
Let $T(\cdot)$ be a bounded analytic semigroup on $X$ generated by $-A$.
Then the system \eqref{eq:control-system} is called 
{\em $L^p$--wellposed of type $\al$}, if 
$C \in B(X_k, Y)$ and $B \in B(U, X_{-k})$ are $L^p$--admissible of type $\al$,  
and $\cF^\al: L^p_\al(U) \to L^p_\al(Y)$, $\cF^\al u = CT_{-k}(\cdot)B
\ast u$ is bounded.
\end{definition}
If we have finite--time $L^p$--admissibility of type $\al$ of $B$ and
$C$ on $I := [0,\tau]$ and if $\cF_\tau^\al: L_\al^p(I, U) \to
L_\al^p(I, Y)$ is bounded, we call the system {\em finite--time
$L^p$--wellposed of type $\al$}. By resorting to $(e^{-\om t}T(t))$ we
may then assume $0\in\varrho(A)$.  The next lemma shows that the only
possible singularity of the convolution kernel of $\cF_\tau$ is at
$t=0$, whence the notion of finite--time $L^p$--wellposedness of type
$\al$ does not depend on $\tau>0$.
Notice that by analyticity of the semigroup $T(\cdot)$, we have
$T_{-k}(t) X_{-k} \subseteq X_{k}$ for $t>0$. Therefore the convolution
kernel $CT_{-k}(\cdot)B$ of $\cF^\al$ is a well defined bounded operator 
from $U$ to $Y$ (see next lemma for norm estimates).
From Theorems \ref{thm:obs-Lp-zul-typ-alpha},
\ref{thm:control-Lp-zul-typ-alpha} and
\ref{thm:beschr-bedingung-und-interpol} we know that
$L^p$--admissibility of type $\al$ of $B$ and $C$

yields $(X, \dot{X}_{k})_{\theta,1} \emb \cD(C)$
and $\cR(B) \subseteq (X, \dot{X}_{-k})_{\si,\infty}$ for
$\theta= (\al+\sfrac{1}p)/k$ and
$\si=(\sfrac{1}{p'}-\al)/k$. Notice that
$k(\si+\theta)=1$. The next lemma is well known to specialists. 
We use it for $q=1$ and $r=\infty$ and provide a proof in
Section~\ref{sec:proofs}. 
\begin{lemma}\label{lem:Abbildungsverhalten-von-analytischen-HG}
  Let $X$ be a Banach space and $T(\cdot)$ be a bounded analytic
  semigroup on $X$. Let $k\in\NN$ and $\si, \theta \in(0,1)$ such that
  $k(\si+\theta)=1$. Let $q,r\in [1,\infty]$  
  and $Z := (X, \dot{X}_k)_{\theta,q}$ and 
  $W := (\dot{X}_{-k},X)_{1-\si,r}$. 
  Then, there exists a constant $M>0$ such
  that  $\norm{ T(t) }_{W \to Z} \le M  / t$ for all $t>0$.
\end{lemma}
To sum up the above considerations: whenever $C$ and $B$ are
$L^p$--admissible of type $\al$, we have 
$\norm{C T(t) B}_{U \to Y} \le M/t$ for $t>0$ for some
constant $M>0$. A corresponding condition in the following theorem
seems thus very natural.
\begin{theorem}\label{thm:PruessSimonett}
Let $p\in(1,\infty)$, $\al\in(-\sfrac{1}p,\sfrac{1}{p'})$ and
$k\in\NN$. Let $X,U,Y$ be Banach spaces and let $T(\cdot)$ be a
bounded analytic semigroup on $X$, $C \in B(X_k, Y)$ and 
$B \in B(U, X_{-k})$ such that $\norm{ CT(t) B}_{U\to Y} \le M/t$. Then 
$\cF^\al := CT(\cdot)B \ast$ is bounded $L^p_\al(U) \to L^p_\al(Y)$ 
if and only if $\cF = \cF^0$ is bounded $L^p(U) \to L^p(Y)$.
\end{theorem}
The above theorem is basically a reformulation of 
\cite[Thm. 2.4]{PruessSimonett}. However, our proof
allows also negative values of $\al$. It relies on the following
\begin{proposition}\label{prop:PruessSimonett}
Let $p\in(1,\infty)$ and $\al<1-\sfrac{1}p$.
Let $U,Y$ be Banach spaces and suppose that $K \in C( \RR_+, B(U,Y))$
satisfies $\norm{K(t)} \le M/t$ for some $M>0$. Let
\[
(Tf)(t) := \int_0^t K(t{-}s) \bigl[ \bigl(\tfrac{t}s\bigr)^\al-1\bigr]\,f(s)\,ds,
\qquad f\in L^p(\RR_+, X).
\]
Then $T \in B( L^p(\RR_+,X), L^p(\RR_+,Y))$ with norm bound $cM$ where
$c=c(p,\al)$.
\end{proposition}
We shortly discuss boundedness of $\cF$ for $\al=0$.
Suppose that $C:Z\to Y$ and $B:U\to W$ are bounded where $Z$ and $W$ are Banach spaces
satisfying $X_1\subseteq Z\subseteq X$ and $X\subseteq W\subseteq \Xm$. Suppose
further that the restriction $A_W$ of $\Am$ to $W$ is sectorial with 
$\cD(A_W)=Z$ (equivalent norms). Then $CT(\cdot)B*:L^p(U)\to L^p(Y)$ is bounded
if $T(\cdot)*:L^p(W)\to L^p(Z)$ is bounded. But it is well--known that the 
latter is equivalent to $A_W$ having the property of \emph{maximal 
$L^p$--regularity} (we refer to
\cite{Amann:Quasilin-Buch,DaPratoGrisvard, Dore:LectureNotes,
      KunstmannWeis:Levico, Weis:FM}
for this relation, the problem of maximal regularity, 
characterisation results and further references on the subject). Thus we 
have obtained
\begin{theorem}
Let $p\in(1,\infty)$ and $W$ be a Banach space $X\subseteq W\subseteq
\Xm$ such that $A_W$ has maximal $L^p$--regularity. Let $C:Z\to Y$ and
$B:U\to W$ be bounded where $Z$ denotes $\cD(A_W)$ equipped with the
graph norm. Then $\cF^\al$ is bounded $L^p_\al(U)\to L^p_\al(Y)$ for
any $\al\in(-\sfrac{1}p,\sfrac{1}{p'})$.
\end{theorem}
Our results can e.g. be used to establish existence and uniqueness of
solutions to some nonlinear systems in feedback form.
Let $p\in (1,\infty)$, let $X$, $U$, $Y$ be Banach spaces and let $x_0\in X$. 
Let $F:X\to B(Y,U)$ be a Lipschitz-continuous operator-valued function, let
$A$, $B$, $C$ be linear operators such that $-A_0:=-(A-BF(x_0)C)$ generates an 
analytic semigroup $T(\cdot)$ in $X$ and such that $B$ and $C$ are
finite--time $L^p$--admissible of type $\al$ for $A_0$. We consider
the closed loop system
\begin{equation}\label{eq:NLP}
\left\{
   \begin{array}{lcl}
     x'(t) + A x(t) &=& B u(t),\\
      x(0)&=&x_0,\\
     y(t)&=&C x(t), \\
     u(t)&=&F(x(t))\, y(t),
   \end{array}
\right.
\end{equation}
on $[0,\tau]$ which we rewrite as
\begin{equation}\label{eq:NLP2}
\left\{
  \begin{array}{lcl}
   x'(t) + A_0 x(t) &=& B \bigl(F(x(t))-F(x_0)\bigr)\, Cx(t) \qquad
   t\in[0,\tau],\\
   x(0) &=&x_0.
  \end{array}\right.
\end{equation}
We are interested in mild solutions of \eqref{eq:NLP2}, i.e. solutions of 
\begin{equation}\label{eq:int-eq}
  x = T(\cdot)x_0+T* B\bigl( F(x)-F(x_0)\bigr)\,Cx\quad\mbox{ on } [0,\tau].
\end{equation}
The following result shall also be proved in Section \ref{sec:proofs}.
\begin{theorem}\label{thm:nonlin}
Assume in  addition to the preceding assumptions
that $Z\emb X$ is a Banach space such that the system corresponding to
$(A_0,B,\text{Id}_Z)$ is finite--time $L^p$--wellposed  of
type $\al$ and that $C\in B(Z,Y)$. Then there exists $\tau=\tau(x_0)>0$
such that \eqref{eq:NLP} has a unique mild solution $x\in
C([0,\tau],X)\cap L^p_\al([0,\tau],Z)$.
\end{theorem}

\begin{remark}
The reason for introducing the space $Z$ here is that $x\mapsto Cx$ may not
induce a closed operator from $C([0,\tau],X)$ into
$L^p_\al([0,\tau],Y)$. If $C$ is closed as an operator from $X$ to
$Y$, then we can replace $C([0,\tau],X)\cap  L^p_\al([0,\tau],Z)$ by
$\{x\in C([0,\tau],X): x(t)\in D(C)\ \mbox{a.e.}\ , Cx\in
L^p_\al([0,\tau],Y) \}$ in the assertion, and the assumption on the 
space $Z$ is not needed.
\end{remark}

\begin{remark}\label{rem:motivation}
Our results on $L^p$--admissibility of type $\al$ and 
$L^p$--wellposedness of type $\al$ give more flexibility in the choice 
of Banach spaces $X$, $U$, and $Y$ for the modelling of a given problem. 
This is important for the study of nonlinear systems via fixed point 
arguments, e.g. via Theorem~\ref{thm:nonlin}, where an appropriate choice 
of $\al$ allows to choose state spaces $X$ as function spaces with little 
regularity (cf. also \cite[Rem. 3.3(b)]{PruessSimonett}).
For the example of the controlled heat equation that we study in
Section~\ref{sec:appl} we refer to Remarks \ref{rem:case-q-gleich-2}
and \ref{rem:case-allgemeines-q} where we describe how to obtain state
spaces with negative smoothness index by suitable choices of $\al$. In 
Example~\ref{ex:nonlin} we give an application of Theorem~\ref{thm:nonlin}
in a nonlinear feedback situation.
We mention in this context that Besov spaces of negative
order have become relevant as spaces for initial values in the study of
other nonlinear partial differential problems, e.g. Navier-Stokes
equations (cf. \cite{Cannone:Buch}). 
\end{remark}
In the next section we provide examples and applications of our
results. In Section \ref{sec:proofs} we shall give proofs of the
results presented so far.

\bigskip
\section{Example: A controlled heat equation}\label{sec:appl}
\medskip

In this section we illustrate our results with a controlled heat equation. 
In \cite{BGSW}, the problem has been studied in the state space 
$X=L^2(\Om)$ and for $\al=0$. Below we discuss Lebesgue and Besov spaces.
Let $\Omega\subset\RR^n$ be a bounded domain with boundary 
$\partial\Omega\in C^\infty$. Denote the outer normal
unit vector on $\partial\Omega$ by $\nu:\partial\Omega\to\RR^n$. 
We consider the following problem
\begin{equation}\label{eq:control-heat}
\left\{
\begin{array}{lcl}
 x'(t)-\Delta x(t) & = & 0,\quad (t>0)\\
 \frac{\partial x(t)}{\partial\nu}\bigr|_{\partial\Omega} & = &
  u(t),\quad (t>0)\\
 x(0) & = & x_0\\
 y(t) & = & x(t)\bigr|_{\partial\Omega},\quad (t>0),
\end{array}
\right.
\end{equation}
where $x(\cdot)$ takes values in some function space $X$. The functions 
$u(\cdot)$ and  $y(\cdot)$ take values in function spaces on the
boundary. For the modelling we follow closely \cite{BGSW} and 
set $A := -\Delta$ with Neumann boundary condition. 
In the state spaces we shall consider below, the operator $A$ is sectorial
of type $0$, but not injective. 
We aim for finite--time $L^c$--admissibility of type $\gamma$
for observation operators  and finite--time
$L^b$--admissibility of type $\beta$ for control operators, where
$b,c \in[1,\infty]$. Sufficient for this is
$L^{b/c}$--admissibility of type $\beta/\gamma$ on $\RR_+$ for
$\Id+A$.  In order to use our characterisations of admissibility,  we
assume
\[
\beta \in \bigl( 
                  -\sfrac{1}{b}, \sfrac{1}{b'} \bigr)
\quad\text{and}\quad
\ga   \in \bigl(-\sfrac{1}{c},  
                                 \sfrac{1}{c'}\bigr).
\]
First, we start with
\bigskip
\subsection*{The $L^q$--case}
Consider $X := L^q(\Om)$, $1<q<\infty$. Due to the smoothness of
$\partial\Omega$ we then have $\cD(A)=\{x\in W^{2}_{q}(\Omega):
\frac{\partial x}{\partial\nu} \bigr|_{\partial\Omega}=0\}$ and
$-A$ generates a bounded analytic semigroup.  
For to ensure \LpAbsch{c}--estimates, we have to impose
\begin{equation}
  \label{eq:Lp-Absch-Bedingung}
  X \emb
  \bigl( \dot{X}_{-1}(\Id+A), \dot{X}_{1}(\Id+A) \bigr)_{\einhalb,c} = 
  \bigl( X_{-1}(A), X_1(A) \bigr)_{\einhalb,c}
= \bigl( X_{-\delta}(A), X_\delta(A) \bigr)_{\einhalb,c}
\end{equation}
where the last equality holds for any $\delta>0$ by reiteration. In
the case $X=L^q(\Om)$, we have $X_\delta=H^{2\delta}_q(\Om)$ for small
$\delta>0$ (cf. \cite{Seeley:interpolation}) and
$X_{-\delta}=H^{-2\delta}_q(\Om)$ by dualisation. Therefore,
\LpAbsch{c}--estimates for $\Id+A$ on $L^q(\Om)$ are equivalent to
the continuous embedding $L^q(\Om) \emb B^0_{q,c}(\Om)$ 
(cf. \cite[Thm 2.4.1, 4.3.1]{Triebel:interpolation}). 
We use the following lemma which shall be proved in Section~\ref{sec:proofs}.
\begin{lemma}\label{lem:kaesekaestchen}
Let $\Omega\subset\RR^n$ be a bounded domain with boundary 
$\partial\Omega\in C^\infty$. If $p,q\in[1,\infty]$, we have $L^q(\Om)
\emb B^0_{q,p}(\Om)$ if and only if $p \ge \max(2,q)$.
\end{lemma}
For an application of Theorem~\ref{thm:obs-Lp-zul-typ-alpha} for 
$L^c$--admissibility of $C$ in $X=L^q(\Om)$ we thus need $c \ge \max(2,q)$.
For an application of Theorem~\ref{thm:control-Lp-zul-typ-alpha} for 
$L^b$--admissibility of $B$ in $X=L^q(\Om)$ in case $\beta=0$, we also need 
\LpAbsch{b'}--estimates for $(I{+}A)'$ on $X'$. By the arguments above
these are equivalent to $L^{q'}(\Om) \emb B^0_{q',b'}(\Om)$ which means 
by Lemma~\ref{lem:kaesekaestchen} that, for $\beta=0$, we have to suppose 
$b \le \min(q,2)$.
Obviously, if $\beta=0$, then application of both Theorems 
\ref{thm:obs-Lp-zul-typ-alpha} and \ref{thm:control-Lp-zul-typ-alpha} 
for $b=c$ would require $b=c=q=2$ and we are back in the Hilbert space
situation. We shall come back to this below.

\medskip
\subsubsection*{Admissibility}
Denoting the Dirichlet trace operator $\gamma_0:x\mapsto
x\bigr|_{\partial\Omega}$ by  $C$ we are looking for a space $Y$ on
$\partial\Omega$ such that $C:\cD(A)\to Y$ is $L^c$--admissible of
type $\ga$ for $\Id+A$.
If $c \ge \max(2,q)$, then $\Id+A$ has \LpAbsch{c}--estimates and we 
know by Theorem~\ref{thm:obs-Lp-zul-typ-alpha} that $L^c$--admissibility
of type $\ga$ of $C:\cD(A)\to Y$ is equivalent to uniform boundedness of
the operators 
\begin{equation*}  
  \la^{1-\ga-\sfrac{1}c} C (\la+\Id+A)^{-1}, \quad \la>0.
\end{equation*}
By Theorem~\ref{thm:beschr-bedingung-und-interpol} this is equivalent to
$C$ having a continuous extension to the Banach space $Z$ where
\begin{eqnarray*}
Z & := & \bigl(X, \dot{X}_{k}(\Id+A)\bigr)_{\ga+\sfrac{1}c,1} \\  
&= & \bigl(X, X_{k}(A) \bigr)_{\ga+\sfrac{1}c,1}   \\
&\overset{(*)}{=}&
\left\{ \begin{array}{ll}
   B^{s}_{q,1}(\Om) 
 & \text{if } s <  1+\sfrac{1}q \\
   \{x\in B^{s}_{q,1}(\Om):\; 
          \tfrac{\partial  x}{d\nu}\bigr|_{\partial\Om} = 0 \} 
 & \text{if } s >  1+\sfrac{1}q,
\end{array}
\right.
\end{eqnarray*}
where $s=2\ga+\sfrac{2}c$. For equality $(\ast)$ we refer to
\cite[Thm 4.3.3]{Triebel:interpolation} or 
\cite[Thm 3.5]{Guidetti:interpolation}.
It is known that $\ga_0$ is bounded from $B^s_{q,1}(\Om)$ to
$B^{s-\sfrac{1}q}_{q,1}(\partial\Om)$ if $s>\sfrac{1}q$
(cf. \cite[Thm. 4.7.1]{Triebel:interpolation}). 
We thus have almost proved the following
\begin{proposition}\label{prop:ex-Lq-obs}
Let $\ga\in\RR$, $c\in[1,\infty]$ and $q\in(1,\infty)$ satisfy 
$c\ge\max(q,2)$ and $2\ga+\sfrac{2}{c}\in(\sfrac{1}{q},1+\sfrac{1}{q})$, and let $Y$ be a Banach 
space and $X=L^q(\Om)$. 
Then $\ga_0:\cD(A)\to Y$ is $L^c$--admissible of type $\gamma$ for $\text{Id}+A$ 
if and only if $B^{2\ga+\sfrac{2}{c}-\sfrac{1}{q}}_{q,1}(\partial\Om)\emb Y$.
\end{proposition}
\begin{proof}
Let $s:=2\ga+\sfrac{2}{c}$. The arguments above show that 
$B^{s-\sfrac{1}{q}}_{q,1}(\partial\Omega)\emb Y$ is sufficient. To prove 
necessity we compose $\ga_0:B^s_{q,1}(\Om)\to Y$ with a continuous 
extension operator
$E_0:B^{s-\sfrac{1}{q}}_{q,1}(\partial\Om) \to B^s_{q,1}(\Om)$ such that 
$\ga_0E_0=Id$
(cf. \cite[Thm. 4.7.1]{Triebel:interpolation}).
\end{proof}
\begin{remark}
The upper bound $2\ga+\sfrac{2}{c}<1+\sfrac{1}{q}$ appears only for simplicity of 
formulation. The calculation of the space $Z$ above indicates how
to proceed in other cases.
\end{remark}

To obtain the representation of the control operator $B$ we follow again 
ideas in \cite{BGSW} and multiply the state equation in 
\eqref{eq:control-heat} with a fixed function  $v\in
C^\infty(\overline{\Omega})$. Then, integrating by parts gives
\[
  \langle x'(t),v \rangle_\Omega + \langle \nabla x(t),\nabla v
  \rangle_\Omega  = \int_{\partial\Omega} u(t) v\, d\sigma,
\]
where $\langle\cdot,\cdot\rangle_\Omega$ denotes the usual duality
pairing on $L^q(\Omega)\times L^{q'}(\Omega)$ and $\sigma$ denotes the
surface measure on $\Gamma:=\partial\Omega$. Denoting extensions of
the usual $L^{q}(\Gamma)\times L^{q'}(\Gamma)$-duality by
$\langle\cdot,\cdot\rangle_\Gamma$ we thus have  
\[
 \int_{\partial\Omega} u(t) v\,d\sigma=\langle u(t),\gamma_0
  v\rangle_\Gamma,
\]
which means that $B=\gamma_0'$ if we identify $X_{-1}(A)$ with  
the dual space of $(X')_1(A')$. 
Notice that $A'=-\Delta$ with Neumann boundary conditions in
$X'=L^{q'}(\Omega)$. 
We are interested in $L^b$--admissibility of type $\beta$ for 
$\ga_0':U\to X_{-1}$ and $\text{Id}+A$ in $X=L^q(\Om)$.
Assuming $b\le\min(q,2)$ in case $\beta=0$ we may use 
Theorem~\ref{thm:control-Lp-zul-typ-alpha} and only have to check
boundedness of $\ga_0':U\to W$ where
the extrapolation space $W$ is given by 
\[
 W  =  (\dot{X}_{-1}(\text{Id}+A),X)_{\beta+\sfrac{1}b,\infty} 
    =  (X_{-1}, X)_{\beta+\sfrac{1}b,\infty} 
    =  \bigl( \bigl( (X')_1,X' \bigr)_{\beta+\sfrac{1}b,1} \bigr)'
    =  \bigl( \bigl( X',(X')_1 \bigr)_{1-(\beta+\sfrac{1}b),1} \bigr)'
\]
(cf. \cite[Sect.3.7]{BerghLoefstroem}).
As we have seen above, we have 
\[
 (X',(X')_{1})_{1-(\beta+\sfrac{1}b),1}=B^{2-2\beta-\sfrac{2}{b}}_{q',1}(\Omega)
\]
for $2-2\beta-\sfrac{2}{b}<1+\sfrac{1}{q'}=2-\sfrac{1}{q}$. 
Now $\ga_0:B^{2-2\beta-\sfrac{2}{b}}_{q',1}(\Omega)
         \to B^{2-2\beta-\sfrac{2}{b}-\sfrac{1}{q'}}_{q',1}(\partial\Omega)$
is continuous for $2-2\beta-\sfrac{2}{b}>\sfrac{1}{q'}=1-\sfrac{1}{q}$ 
and $\partial\Omega$ is without boundary, hence
$(B^s_{q',1}(\partial\Omega))'=B^{-s}_{q,\infty}(\partial\Omega)$.
Thus we have proved one implication in the following 
\begin{proposition}\label{prop:ex-Lq-control}
Let $\beta\in\RR$, $b\in[1,\infty]$ and $q\in(1,\infty)$ satisfy 
$\sfrac{2}{b}+2\beta\in(\sfrac{1}{q},1+\sfrac{1}{q})$. If $\beta=0$ assume in addition
$b\le\min(q,2)$. Let $U$ be a Banach space and $X=L^q(\Om)$. 
Then $\ga_0':U\to X_{-1}$ is $L^b$--admissible of type $\beta$ for $\text{Id}+A$ 
if and only if $U\emb B^{\sfrac{2}{b}+2\beta-1-\sfrac{1}{q}}_{q,\infty}(\partial\Om)$.
\end{proposition}
For the remaining implication we make use of $E_0'$ where $E_0$ is the
extension map from the proof of Proposition~\ref{prop:ex-Lq-obs}.

\medskip
\subsubsection*{Wellposedness of the full system}
We study $L^p$--wellposedness of type $\al$ in finite time for the
system \eqref{eq:control-heat} in the state space $X=L^q(\Om)$ where 
$\al\in\RR$, $p\in[1,\infty]$ and $q\in(1,\infty)$.
Hence we have $c=b=p$ and $\ga=\beta=\al$ in the admissibility
situations above. Again, it is sufficient to study infinite--time
$L^p$-admissibility of type $\al$ for $Id+A$.
The assumptions of Propositions \ref{prop:ex-Lq-obs}
and \ref{prop:ex-Lq-control} lead, for $\al\neq0$, to the restrictions
$p\ge\max(q,2)$ and $2\al+\sfrac{2}p \in(\sfrac{1}q,1+\sfrac{1}q)$,
and for $\al=0$ to the restriction $p=q=2$.

For $p=q=2$ we thus can state the following result extending 
the case $\al=0$, $r=2$, considered in \cite{BGSW}.
\begin{proposition}\label{prop:ex-Lq-wellp}
Let $|\al|<\sfrac{1}4$, $r\in[1,\infty]$ and $X=L^2(\Om)$. 
Let $Y$ and $U$ be Banach spaces on $\partial\Om$ such that 
$B^{2\al+\einhalb}_{2,r}(\partial\Om)\emb Y$ and
$U\emb B^{2\al-\einhalb}_{2,r}(\partial\Om)$. Then the system
\eqref{eq:control-heat} is $L^2$--wellposed of type $\al$.
\end{proposition}
\begin{proof}
By Propositions \ref{prop:ex-Lq-obs} and \ref{prop:ex-Lq-control},
$C$ and $B$ are $L^p$--admissible of type $\al$ 
for $\text{Id}+A$. Notice that the semigroup associated
to $\text{Id}+A$ is $S(t) := e^{-t} T(t)$. 
Since $B=\ga_0:B^{2\al+\einhalb}_{2,r}(\partial\Om)\to
\widetilde{W}:=(\Xm,X)_{\beta+\sfrac{1}{b},r}
=(\Xm,X)_{\al+\einhalb,r}$ (recall $\beta=\al$ and $b=p=2$) and
$C=\ga_0:\widetilde{Z}\to B^{2\al+\einhalb}_{2,r}(\partial\Om)$ where
$\widetilde{Z}:=B^{2\al+1}_{2,r}(\Om)=(X,X_1)_{\al+\einhalb,r}$ it 
suffices to prove that 
$S(\cdot)\ast: L_\al^p(\RR_+, \widetilde{W}) \to L_\al^p(\RR_+,\widetilde{Z})$.
Theorem~\ref{thm:PruessSimonett} applies by $|\al|<\sfrac{1}4<\einhalb$.
Now notice that by $\widetilde{Z}=\cD(A_{\widetilde{W}})$. Hence we are left 
to check that $A_{\widetilde{W}}$ has maximal $L^p$--regularity in 
$\widetilde{W}$ which holds by \cite[Thm. 4.7]{DaPratoGrisvard}.
\end{proof} 
With essentially the same arguments we can study 
problem \eqref{eq:control-heat} in state spaces $X=H^\delta_q(\Om)$ where we 
restrict to $\delta\in(-\sfrac{1}{q'},\sfrac{1}q)$ and are thus not bothered 
by additional boundary conditions. 
Here we have $X_1=\cD(A)=\{x\in H^{2+\delta}_q(\Om):
\frac{\partial x}{\partial\nu}|_{\partial\Om}=0\}$, and a repetition of the
arguments above yields for the corresponding interpolation spaces 
$Z=B^{\delta+2(\al+\sfrac{1}p)}_{q,1}(\Om)$ 
if $\delta+2(\al+\sfrac{1}p)<1+\sfrac{1}q$
and $W=(B^{2-\delta-2(\al+\sfrac{1}p)}_{q',1}(\Om))'$ 
if $2-\delta-2(\al+\sfrac{1}p)<1+\sfrac{1}{q'}$, 
i.e. if $\sfrac{1}q<\delta+2(\al+\sfrac{1}p)$.
Thus we obtain
\begin{proposition}\label{prop:ex-Hdel2}
Let $q\in(1,\infty)$ and $p\in[1,\infty]$ satisfy $p\ge\max(2,q)$.
Let $\delta \in (-\sfrac{1}{q'},\sfrac{1}{q})$, 
and $\delta+2(\al+\sfrac{1}p)\in(\sfrac{1}q,1+\sfrac{1}q)$ where $\al\neq0$. 
Let $X=H^\delta_q(\Om)$ and let $Y$ and $U$ be Banach spaces on $\partial\Om$. 
\begin{enumerate}
\item
The operator $C=\ga_0: \cD(A)\to Y$ is $L^p$--admissible of type $\al$
for $\Id{+}A$ if and only if\\$B^{\delta+2(\al+\sfrac{1}p)-\sfrac{1}q}_{q,1}(\partial\Om)\emb Y$.
\item
The operator $B=\ga_0': U\to \Xm$ is $L^p$--admissible of type $\al$
for $\Id{+}A$ if and only if\\$U\emb B^{\delta+2(\al+\sfrac{1}p)-\sfrac{1}q-1}_{q,\infty}(\partial\Om)$. 
\item
Let $r\in[1,\infty]$ and 
$B^{\delta+2(\al+\sfrac{1}p)-\sfrac{1}q}_{q,r}(\partial\Om)\emb Y$ 
and $U\emb B^{\delta+2(\al+\sfrac{1}p)-\sfrac{1}q-1}_{q,r}(\partial\Om)$.\\ 
Then the system \eqref{eq:control-heat} is finite--time $L^p$--wellposed of type $\al$.
\end{enumerate}
Observe that the assumptions imply $\al+\sfrac{1}p \in (0,1)$.
\end{proposition}
We see that $\delta$ may be chosen arbitrarily close to $-\sfrac{1}{q'}$ by taking
$p$ large and adjusting $\al$. Moreover, one still has free choice of $r\in[1,\infty]$.

\bigskip
\subsection*{The Besov space case}
We continue the study of \eqref{eq:control-heat}, but now we take
as state space the Besov space $X := B_{q,v}^0(\Omega)$ where 
$1<q,v<\infty$ are fixed. 
Again we put $A:=-\Delta$ with homogeneous Neumann boundary 
conditions, i.e., $A$ has domain 
$\cD(A)=\{x\in B^{2}_{q,v}(\Omega): \frac{\partial x}{\partial\nu}
\bigr|_{\partial\Omega}=0\}$. As in the $L^q$--case, $-A$ generates a
bounded analytic semigroup  in $X$.

\subsubsection*{Admissibility}
Still denoting by $\ga_0$ the Dirichlet trace on $\partial\Om$, 
we study $L^c$--admissibility of type $\ga$ for $C=\ga_0$ and 
$L^b$--admissibility of type $\beta$ for $B=\ga_0'$ where 
$\beta,\ga\in\RR$ and $b,c\in[1,\infty]$.
First we note that 
$\bigl( X_{-1}, X_{1} \bigr)_{\einhalb,c} =
B_{q,c}^{0}(\Om)$, whence  $\Id{+}A$ has \LpAbsch{c}--estimates on
$X$ if and only if $B_{q,v}^0(\Om) \emb B_{q,c}^0(\Om)$ which is 
equivalent to $c\ge v$. Notice that this condition does not depend 
on $q$. 
If $c\ge v$, then the characterising Theorem~\ref{thm:obs-Lp-zul-typ-alpha}
applies, and $L^c$--admissibility of type $\ga$ 
of the observation operator $C=\ga_0$ may be checked by verifying 
the boundedness condition on the set $W_C$. 
Again, we use Theorem~\ref{thm:beschr-bedingung-und-interpol} and find that,
for $2\ga+\sfrac{2}{c}<1+\sfrac{1}{q}$, $W_C$ is bounded if and only if $\ga_0$ has 
a continuous extension to the very same Banach space 
\[
 Z = B^{2\ga+\sfrac{2}{c}}_{q,1}(\Om)
\]
we calculated above. Hence the proof of the following can be done as 
in the case $X=L^q(\Om)$.
\begin{proposition}\label{prop:ex-Bqr-obs}
Let $\ga\in\RR$, $c\in[1,\infty]$ and $q,v\in(1,\infty)$ satisfy 
$c\ge v$ and $2\ga+\sfrac{2}{c}\in(\sfrac{1}{q},1+\sfrac{1}{q})$. 
Let $Y$ be a Banach space and $X=B_{q,v}^0(\Om)$. 
Then $\ga_0:\cD(A)\to Y$ is $L^c$--admissible of type $\ga$ for $\text{Id}+A$ 
if and only if $B^{2\ga+\sfrac{2}{c}-\sfrac{1}{q}}_{q,1}(\partial\Om)\emb Y$.
\end{proposition}
We turn to $L^b$--admissibility of type $\beta$ for $B=\ga_0'$.
For an application of Theorem~\ref{thm:control-Lp-zul-typ-alpha} in
$X=B_{q,v}^0(\Om)$ in case $\beta=0$ we need \LpAbsch{b'}--estimates
for $(\Id{+}A)'$, which  are, by the argument above, equivalent to
$b\le v$. This yields
\begin{proposition}\label{prop:ex-Bqr-control}
Let $\beta\in\RR$, $b\in[1,\infty]$ and $q,v\in(1,\infty)$ satisfy 
$\sfrac{2}{b}+2\beta\in(\sfrac{1}{q},1+\sfrac{1}{q})$. 
Assume additionally $b\le v$ if $\beta=0$.
Let $U$ be a Banach space and $X=B_{q,v}^0(\Om)$. 
Then $\ga_0':U\to X_{-1}$ is $L^b$--admissible of type $\beta$ for $\text{Id}+A$ 
if and only if $U\emb B^{\sfrac{2}{b}+2\beta-1-\sfrac{1}{q}}_{q,\infty}(\partial\Om)$.
\end{proposition}
\medskip
\subsubsection*{Wellposedness of the full system}
We study $L^p$--wellposedness of type $\al$ for the system 
\eqref{eq:control-heat} where $\al\in\RR$, $p\in[1,\infty]$ 
and the state space $X=B_{q,v}^0(\Om)$ with $q,v\in(1,\infty)$.
Again we have $c=b=p$ and $\ga=\beta=\al$ in the admissibility
situations above. The assumptions of Propositions \ref{prop:ex-Bqr-obs}
and \ref{prop:ex-Bqr-control} now lead, for $\al\neq0$, to the restrictions
$p\ge v$ and $2\al+\sfrac{2}{p} \in (\sfrac{1}{q},1+\sfrac{1}{q})$.
In case $\al=0$ we are led to $p=v$. 
For example, if $p=q$, then we are led to 
$2\al \in (-\sfrac{1}{q},\sfrac{1}{q'})$
which, for $q=2$, corresponds to the condition $|\al|<\sfrac{1}4$
in Proposition~\ref{prop:ex-Lq-wellp} (recall $L^2(\Om)=B_{2,2}^0(\Om)$).
By the arguments used above we hence obtain
\begin{proposition}\label{prop:ex-Bqr-wellp}
Let $\al\in\RR$, $p\in(1,\infty)$ and $q\in(1,\infty)$ satisfy 
$2\al+\sfrac{2}{p} \in(\sfrac{1}{q}, 1+\sfrac{1}{q})$.
Let $X=B_{q,p}^0(\Om)$ and let $Y$ and $U$ be Banach spaces 
on $\partial\Om$ such that 
$B^{2\al+\sfrac{2}{p}-\sfrac{1}{q}}_{q,p}(\partial\Om)\emb Y$ and
$U\emb B^{\sfrac{2}{p}+2\al-1-\sfrac{1}{q}}_{q,p}(\partial\Om)$. 
Then the system \eqref{eq:control-heat} is finite--time
$L^p$--wellposed of type $\al$.
\end{proposition}
\begin{proof}
The proof is similar to the case $X=L^q(\Om)$. We check here that 
Theorem~\ref{thm:PruessSimonett} applies, i.e., 
that $\al\in(-\sfrac{1}p,\sfrac{1}{p'})$. The condition on 
$\al$ in the assumption may be rephrased as 
$\al\in(\sfrac{1}{2q}-\sfrac{1}{p},\einhalb+\sfrac{1}{2q}-\sfrac{1}{p})$.
In particular $\al>\sfrac{1}{2q}-\sfrac{1}p\ge-\sfrac{1}p$ and
$\al<\einhalb+\sfrac{1}{2q}-\sfrac{1}{p}\le 1-\sfrac{1}p$.
\end{proof}
We also give an analogue of Proposition~\ref{prop:ex-Hdel2} for Besov spaces, e.g.
for state spaces $X=B^\delta_{q,p}(\Om)$ where $q,p\in(1,\infty)$ and 
we restrict to $\delta\in(-\sfrac{1}{q'},\sfrac{1}q)$ for the same reasons 
as before.
\begin{proposition}\label{prop:ex-Besov-del}
Let $q,p\in(1,\infty)$, $\delta\in(-\sfrac{1}{q'},\sfrac{1}q)$,
$2\al+\sfrac{2}p+\delta\in(\sfrac{1}q,1+\sfrac{1}q)$,
and $X=B^\delta_{q,p}(\Om)$. 
Let $Y$ and $U$ be Banach spaces on $\partial\Om$. 
\begin{enumerate}
\item
The operator $C=\ga_0: \cD(A)\to Y$ is $L^p$--admissible of type $\al$
for $\Id{+}A$ if and only if\\ $B^{2\al+\sfrac{2}p+\delta-\sfrac{1}q}_{q,1}(\partial\Om)\emb Y$.
\item
The operator $B=\ga_0': U\to \Xm$ is $L^p$--admissible of type $-\al$
for $\Id{+}A$ if and only if\\ 
$U\emb B^{2\al+\sfrac{2}p+\delta-1-\sfrac{1}q}_{q,\infty}(\partial\Om)$. 
\item
Let $r\in[1,\infty]$ and 
$B^{2\al+\sfrac{2}p+\delta-\sfrac{1}q}_{2,r}(\partial\Om)\emb Y$ 
and $U\emb B^{2\al+\sfrac{2}p+\delta-1-\sfrac{1}q}_{2,r}(\partial\Om)$.\\
Then the system \eqref{eq:control-heat} is finite--time $L^p$--wellposed of type $\al$.
\end{enumerate}
\end{proposition}
\subsection*{Discussion}
We discuss our results for the system \eqref{eq:control-heat} by
starting from $Y:=B^s_{q,r}(\partial\Om)$ and
$U:=B^{s-1}_{q,r}(\partial\Om)$ where $s\in(0,1)$, $q\in(1,\infty)$
and $r\in[1,\infty]$ are fixed. We now look for $\al$, $p$ and a state
space $X$ such that \eqref{eq:control-heat} is $L^p$--wellposed of type $\al$.
This should be compared with the situation studied in \cite{BGSW} where
$p=q=r=2$, $\al=0$, $s=\einhalb$, and $X=L^2(\Om)$.
\begin{remark}\label{rem:case-q-gleich-2}
In case $q=2$, we may take $p=2$ and $X=H^\delta_2(\Om)$ where the
restrictions may be  read off of Proposition~\ref{prop:ex-Hdel2}: 
$2\al+\delta+\einhalb=s$ and $|\delta|<\einhalb$.
Thus we see that a suitable choice of $\al$ always allows to have
$\delta$ arbitrarily close to $-\einhalb$.
In particular this applies to the ``classical'' case $s=\einhalb$ where
the restriction $\al=0$ forces $\delta=0$ and $X=L^2(\Om)$.
\end{remark}

\begin{remark}\label{rem:case-allgemeines-q}
For general $q\in(1,\infty)$, we may apply Proposition~\ref{prop:ex-Besov-del}
and take $X=B^\delta_{q,p}(\Om)$  under the restrictions
$s=2\al+\sfrac{2}p+\delta-\sfrac{1}{q}$ and
$\delta\in(-\sfrac{1}{q'},\sfrac{1}q)$. We see that, taking $p$ arbitrarily
large and $\delta$ arbitrarily close to $-\sfrac{1}{q'}$, we still
obtain finite--time $L^p$--wellposedness of type $\al$ for the system
with state space $X=B^\delta_{q,p}(\Om)$ by adjusting $\al$.
Observe that the state space $X=H^\delta_q(\Om)$ would have required
the additional restriction $p\ge\max(2,q)$ for $\al\neq0$.
\end{remark}

\begin{remark}\label{rem:ex-nonlin}
We let $s=\einhalb$ and $U := H^{-\einhalb}_2(\partial\Omega)$, 
$Y := H^{\einhalb}_2(\partial\Omega)$, i.e. we take $q=r=2$ in the situation 
discussed above. Then, for $\delta\in(-\einhalb,\einhalb)$ and 
$X=B^\delta_{2,p}(\Om)$, we have that \eqref{eq:control-heat} is
finite--time $L^p$--wellposed 
of type $\al$ if $\delta=1-2(\al+\sfrac{1}p)$. The restrictions are 
$p\in(1,\infty)$ and $\al+\sfrac{1}p\in(\sfrac{1}4,\sfrac{3}4)$.
This means in other words that, for $\eps\in(0,1)$, $p\in[2,\infty)$,
and $X=B^{-\einhalb+\eps}_{2,p}(\Om)$, the system \eqref{eq:control-heat} is 
finite--time $L^p$--wellposed of type $\al=\sfrac{3}4-\sfrac{\eps}2-\sfrac{1}p$.
\end{remark}
Now we consider a nonlinear feedback in the setting of the above remark. 
\begin{example}\label{ex:nonlin}
Since $Y=H^\einhalb_2(\partial\Omega)\emb L^{\frac{2n-2}{n-2}}(\partial\Om)$ 
and $L^{2-\sfrac{2}n}(\partial\Om)\emb H^{-\einhalb}_2(\partial\Om)$, H\"older's
inequality yields that any $g\in L^{n-1}(\partial\Om)$ induces a 
bounded multiplication operator $Y\to U, y\mapsto g\cdot y$. We take a smooth
open subset $\Om_0$ with $\overline{\Om_0}\subset\Om$, e.g. a small ball,
and let $\psi(x):=\int_{\Om_0}x(\om)\,d\om=\langle x, 1_{\Om_0}\rangle$
for $x\in X$. Observe that $\psi\in X'$ for all spaces $X$ we mentioned above
(cf. \cite{Triebel:interpolation}).  
Taking a Lipschitz-continuous function $f:\RR\to\RR$, we let 
$F(x)y:=f(\psi(x))\,g\cdot y$. Then $F:X\to B(Y,U)$ is Lipschitz-continuous. 
We interpret $x \mapsto\psi(x)$ as a distributed measurement in $\Om$ 
affecting via $f(\psi(x))$ the intensity of the linear feedback $y
\mapsto g\cdot y$.
For $\eps\in(0,1)$, $p\in[2,\infty)$, $X=B^{-\einhalb+\eps}_{2,p}(\Om)$, 
and $\al=\sfrac{3}4-\sfrac{\eps}2-\sfrac{1}p$ we take $x_0\in X$ satisfying
(for simplicity) $f(\psi(x_0))=0$ and show that Theorem~\ref{thm:nonlin} 
applies. We have $F(x_0)=0$, hence $A_0=A$. We take $Z=H^1_2(\Om)$
whence by an application of Theorem~\ref{thm:obs-Lp-zul-typ-alpha} and
Theorem~\ref{thm:beschr-bedingung-und-interpol}, $\text{Id}_Z$ is
finite--time $L^p$--admissible of type $\al$ (cf. \cite[Remark
2.8.1]{Triebel:interpolation}). 
Since $B=\ga_0': H^{-\einhalb}_2(\partial\Om)\to (H^1_2(\Om))'=:W$ is bounded
and $\Id+A$ has maximal $L^p$-regularity in $W$ we obtain by
Remark~\ref{rem:dualisieren-von-lp-al-zulaessig} and Theorem~\ref{thm:PruessSimonett}
that the system $(A_0,B,\text{Id}_Z)$ is $L^p$--wellposed of type
$\al$ in finite time. Hence we may apply Theorem~\ref{thm:nonlin} and obtain
that the nonlinear system
\begin{equation}
\left\{
   \begin{array}{lcl}
     x'(t) - \Delta x(t) &=& 0, \qquad t\in[0,\tau], \\
     x(0)&=&x_0,\\
     \frac{\partial x(t)}{\partial\nu}
     &=& f(\int\limits_{\Om_0}x(t)\,d\om)\, g \cdot x(t)|_{\partial\Om}, 
   \end{array}
\right.
\end{equation}
has, for some $\tau(x_0)>0$, a unique mild solution $x(\cdot)$ in 
$C([0,\tau],B^{-\einhalb+\eps}_{2,p}(\Om))\cap L^p_\al([0,\tau],H_2^1(\Om))$.
This means that we obtain solutions also for rather rough initial data $x_0$.
\end{example}

In the example above we made the assumption $f(\psi(x_0))=0$ for simplicity.
For initial data $x_0$ with $f(\psi(x_0))\neq 0$ one may resort to perturbation
results in \cite{HaakHaaseKunstmann}.

\bigskip
\section{Proofs}\label{sec:proofs}
\subsection*{Finite--time admissibility of type $\al$}
\begin{proof}[Proof of {Lemma~\ref{lem:alternative-bed-fuer-obs-zul}}]
Assume $C$ to be finite--time $L^p$--admissible of type $\al$ for $A$. For
$b \ge 0$ we then have
\begin{eqnarray*}
      \biggl( \int_{b+\sfrac{\tau}2}^{b+\tau} \norm{ t^\al CT(t) x }^p\,dt \biggr)^{\sfrac{1}p}
&=&   \biggl( \int_{\sfrac{\tau}2}^\tau \norm{ s^\al(1+\sfrac{b}s )^\al CT(s)
       T(b) x }^p\,dt \biggr)^{\sfrac{1}p}\\
&\le& \biggl(\int_0^\tau \norm{ s^\al CT(s) T(b)x }^p\,dt
       \biggr)^{\sfrac{1}p} \, \max((1+\tfrac{2b}\tau)^\al,1)\\
&\le& M_\tau\,\norm{ T(b) x }\,\max((1+\tfrac{2b}\tau)^\al,1).
\end{eqnarray*}
If $T(\cdot)$ is uniformly exponentially stable, i.e. if there are $c, \eps>0$
such that $\norm{ T(t) } \le c\,e^{-\eps t}$, $t\ge 0$, then
we obtain, for  $b=(k{-}1)\sfrac{\tau}2$ and $k\in\NN_0$,
\[
 \biggl( \int_{k\sfrac{\tau}2}^{(k+1)\sfrac{\tau}2} \norm{ t^\al CT(t)
   x }^p\,dt \biggr)^{\sfrac{1}p}
 \le c \,M_\tau\, \max((1+\tfrac{2b}\tau)^\al,1) \, e^{-\eps (k-1)\sfrac{\tau}2}\norm{ x }.
\]
Thus \eqref{eq:def-Lp-al-zul-obs} holds.
A similar reasoning shows that the notion of finite--time
$L^p$--admissibility of type $\al$ for control operators is
independent of $\tau$: let $a = \sfrac{\tau}2$ and $b =
\sfrac{3}2\,\tau$.  Since
\[
     \int_0^{b} T(b{-}s) B u(s)\,ds
\le  T(a) \int_0^{\tau} T(\tau{-}s) B u(s)\,ds
  +  \int_{a}^{\tau} T(\tau{-}s) B u(s{+}a)\,ds,
\]
$L^p$--admissibility of type $\al$ on $[0,\tau]$ gives for $p<\infty$
\begin{eqnarray*}
      \biggnorm{ \int_0^{b} T(b{-}s) B  u(s)\,ds }
&\le& K_\tau \norm{ T(a) } \; \norm{u}_{L^p_\al([0,\tau])}
    + K_\tau \norm{ \eins_{[a, \tau]}(\cdot)
      u(a{+}\cdot) }_{L^p_\al([0,\tau])} \\
&\le& K_\tau \norm{ T(a) } \; \norm{u}_{L^p_\al([0,\tau])}
    + K_\tau \widetilde c_\al \,
      \norm{ u }_{L^p_\al([\tau, b])}\\
&\le& K_\tau \bigl( \widetilde c_\al + \norm{ T(\sfrac{\tau}2) }\bigr) \, 
      \norm{ u }_{L^p_\al([0, b])},
\end{eqnarray*}
where $\widetilde c_\al = \max\bigl((\tfrac12)^\al, (\tfrac23)^\al
\bigr)$ does not depend on $\tau$. In case $p=\infty$, we obtain the
same estimate with $\widetilde c_\al = 1$ by directly regarding the first
line of the above inequalities. Thus, $B$ is $L^p$--admissible of type
$\al$ on $[0,\sfrac{3}2\, \tau]$ as required. An iteration of the
argument shows that
\[
  K_{(\sfrac32)^n \tau} \le K_\tau (\widetilde c_\al)^n \prod_{j=1}^n \bigl(1+
  \widetilde c_\al^{-1} \norm{ T( 3\tau (\sfrac32)^j ) } \bigr).
\]
If $T(\cdot)$ is uniformly exponentially stable and if $\al \ge
0$ or if $p=\infty$, the left hand side of the above inequality is bounded
since  $\sum\limits_{j\ge 1} \exp(-3 \eps \tau (\sfrac32)^j ) < \infty$ and
$|\widetilde c_\al| \le 1$.
\end{proof}
If $\al <0$ and $p<\infty$, equivalence of finite--time and inifinite--time
$L^p$--admissibility of type $\al$ fails in general:
\begin{example}\label{ex:gegenbeispiel-alpha-negativ}
Take $U=X=\CC$, $-\al=\beta>0$ and $T(t)=e^{-\eps t}$. 
Then, for any $\tau>0$, 
\[
    | T*u(\tau) | 
\le \norm{ u }_{L^p_\al(0,\tau)} \norm{ e^{-\eps(\cdot)} }_{L^{p'}_\beta(0,\tau)}
\le c_{p'} \tau^{\beta+\sfrac{1}{p'}} \norm{ u }_{L^p_\al(0,\tau)},
\]
and the identity is finite--time $L^p$--admissible of type $\al$.
On the other hand, letting $u_k := \eins_{[k,k+1]}e^{-\eps(\cdot-k)}$ for $k\in\NN$,
we have $\norm{ u_k }_{L^p_\al} \le k^{-\beta}\norm{ u_k }_{L^p} =c k^{-\beta}$
and for $\delta \in (0,1]$, 
\[
   (T * u_k)(k+\delta) = \int_k^{k+\delta} e^{-\eps(k+\delta-s)}\,e^{-\eps(s-k)}\,ds
   = \delta e^{-\eps \delta}.
\]
Hence the identity is not infinite--time $L^p$--admissible of type $\al$.
\end{example}

\subsection*{Characterisation of $L^p$--admissibility of type $\al$}
In the proof of Theorems \ref{thm:obs-Lp-zul-typ-alpha} and
\ref{thm:control-Lp-zul-typ-alpha} we make use of the following 
lemma.
\begin{lemma}[{\cite[Lem. 4.1]{HaakLeMerdy}}]\label{lem:technical}
Let $\sigma\in (0,\pi)$, let $\varphi\in
H^{\infty}_{0}(\Sec{\sigma})$, and let $m\geq 1$ be an integer.
There exist a function $f \in H^{\infty}_{0}(\Sec{\sigma})$ and a
constant $a\in\CC$ such that
\begin{equation}\label{eq:darstellung}
\varphi(z) = z^m\, f^{(m)}(z)\, +\, a \,\frac{z^m}{(1+z)^{m+1}}
\,,\qquad z\in\Sec{\sigma}.
\end{equation}
Furthermore, if $\delta,\epsilon\in(0,1)$ are positive numbers
such that
\begin{equation*}
\vert \varphi(z)\vert = O(\vert z\vert^{-\delta})\quad\hbox{at}\
\infty\qquad\hbox{and}\qquad \vert \varphi(z)\vert = O(\vert
z\vert^{\epsilon})\quad\hbox{at}\ 0,
\end{equation*}
then $f$ can be chosen so that we also have $\vert f(z)\vert =
O(\vert z\vert^{-\delta})$ at $\infty$, and $\vert f(z)\vert =
O(\vert z\vert^{\epsilon})$ at $0$.
\end{lemma}

\begin{proof}[Proof of {Theorem~\ref{thm:obs-Lp-zul-typ-alpha}}]
{\it (a)}
Notice that for any $\lambda\in \CC$ with positive
real part and for any $x\in X$, we have
\[
(\lambda+A)^{-k} x = \frac{1}{(k-1)!}\,\int_{0}^{\infty}
t^{k-1} e^{-\lambda t} T(t)x\, dt.
\]
For $x\in X_k=\cD(A^k)$, the integrand $t \mapsto t^{k-1} e^{-\la k}T(t)
x$ belongs to $L^1(\RR_+,X_k)$, and so continuity of $C$ on $X_k$ shows
\[
C (\lambda+A)^{-k} x = \frac{1}{(k-1)!}\, \int_{0}^{\infty}
t^{k-1} e^{-\lambda t} C T(t)x\, dt.
\]
By Hölder's inequality we thus have for $x\in X_k$
\begin{eqnarray*}
       \bignorm{ C(\lambda+A)^{-k} x } 
&\le & \frac{1}{(k-1)!} \, \int_{0}^{\infty} 
       \bignorm{t^{\alpha} CT(t)x } t^{k-1-\alpha} e^{-\Re(\la)\,t} \,dt\\
&\le & \frac{1}{(k-1)!} \biggnorm{ t\mapsto t^\al CT(t)x }_{L^p(\RR_+,Y)}
       \biggl( \int_0^\infty t^{(k-1-\al)p'} e^{-\Re(\la)p'\, t}
       \,dt\biggr)^{\sfrac{1}{p'}}\\
&\overset{s=\Re(\la)p'\, t}{\le} &  \frac{M}{(k-1)!} \norm{ x }_X
       \biggl( (p'\Re(\la))^{-1-(k-1-\al)p'} 
       \Ga\bigl(1{+}(k{-}1{-}\al)p'\bigr) \biggr)^{\sfrac{1}{p'}}\\
& = &  K  \Re(\la)^{-k+\al+\sfrac{1}p} \; \norm{x}_X,
\end{eqnarray*}
where $\Gamma$ is the usual Gamma function and the
number $K$ depends only on $k$, $p$ and the admissibility constant
$M$. By density of $X_k$ in $X$ this shows the first assertion.
\medskip
{\it (b)}
Without loss of generality we may assume $k\ge 2$ since by
sectoriality of $A$, whenever $W_C$ is bounded for some $k\in \NN$, it
is also bounded when $k$ is replaced by $k+1$ (see Remark~\ref{rem:k-egal}).
We make use of the (unbounded) operator $A^{-1}$ that is densely
defined on the range of $A$. We set $F_k(z) := z^{k-1} e^{-z}$. Then for
any $x\in X_{k}$ and any $t>0$, we have
\begin{equation}\label{eq:decomp}
t^{\alpha} C T(t) x = t^{\alpha-k+1} C A^{-k+1} F_k(tA) x.
\end{equation}
For some $\eps\in (0,1)$ that we will precise later on 
consider the decomposition $F_k(z)=\varphi(z)\psi(z)$ where
\begin{equation}\label{eq:decomp2}
\varphi(z) = z^\eps(1+z)^{-1}, \qquad\text{and}\qquad
\psi(z)=z^{k-1-\eps}(1+z)  e^{-z}.
\end{equation}
Note that $\psi\in H_0^\infty(\Sec{\theta})$ for any
$\theta<\frac{\pi}{2}$, whereas $\varphi\in H_0^\infty(\Sec{\sigma})$
for any $\sigma<\pi$. By \eqref{eq:decomp}, we have
\begin{equation*}
    \int_0^\infty \bignorm{ t^\alpha  CT(t) x}_Y^p dt 
\le \int_0^\infty \bignorm{ \underset{=:K(t)}
  {\underbrace{\bigl[t^{\alpha+\sfrac{1}p-k+1} CA^{-k+1}\psi(tA) \bigr]}}
  t^{-\sfrac{1}p}\varphi(tA)x}_X^p \,dt.
\end{equation*} 
Now we will show that the operator family $K(t)$, $t>0$, is uniformly
bounded. Once this is done, the assertion of the theorem follows
immediately from the assumed \LpAbsch{p}--estimate for $A$
(cf. (\ref{eq:Lp-stern-abschaetzungen-und-zul})).
\medskip
We fix $\sigma\in (\om,\pi)$ and apply Lemma~\ref{lem:technical} 
to $\psi$ with $m=k-1$ and $\delta=1-\eps$. Let 
$f\in H_0^\infty(\Sec{\sigma})$
denote the corresponding function satisfying
equation~\ref{eq:darstellung}. Note that according to that equation, $z\mapsto
z^{k-1} f^{(k-1)}(z)$ belongs to $H^{\infty}_{0}(\Sec{\sigma})$. Let
$\theta=\theta_\si$ for some $\theta\in(\omega,\sigma)$ and let $\Ga$
denote the positively orientated boundary of $\Sec{\theta}$.
Then, as in the proof of \cite[Thm. 4.2]{HaakLeMerdy} the following
representation formula for $x$ in the dense subspace $Z :=
\ran{A^{k-1}(I{+}A)^{-k}}$ holds:
\begin{equation}\label{eq:repres}
CA^{-k+1} [z^{k-1} f^{(k-1)}(z)](tA)x = \tfrac{(k-1)!}{2\pi i} \int_\Ga
f(\la) t^{k-1}\,CR(\la,tA)^{k}x \,d\la, \quad t>0.
\end{equation}
For $\la\in \Ga$, by the resolvent equation we have
\begin{equation}
  \label{eq:beschraenktheit-zieht-sich-auf-Gamma-hoch}
  \la^{k-\al-\sfrac{1}{p}} C R(\la,A)^{k} = |\la|^{k-\al-\sfrac{1}{p}} C
(|\la|+A)^{-k} \; \bigl[ 2 \cosh(\pm\sfrac{\theta}2) \la R(\la,A) -I \bigr]^{k},
\end{equation}
and thus $\la^{k-\al-\sfrac{1}p} C R(z,A)^{k}$ is uniformly
bounded by sectoriality of $A$. Now, by the representation
\eqref{eq:repres},
\[
K(t) = \tfrac{(k-1)!}{2\pi i} \int_\Ga t^{\al+\sfrac{1}p} f(\la) \,CR(\la,tA)^{k}x
\,d\la + a t^{\al+\sfrac{1}p} C(I+tA)^{-k}\\
\]
on $Z$. Next we show that for an appropriately chosen $\eps\in(0,1)$,
$K(t)\in B(X,Y)$  and moreover  the operators $K(t)$ are uniformly
bounded for $t>0$. To this end, write
\begin{eqnarray*}
K(t)
& =& \tfrac{(k-1)!}{2\pi i} \int_\Ga f(\la) \la^{1-k+\al+\sfrac{1}p}
\,\left[\left(\tfrac{\la}t\right)^{k-\al-\sfrac{1}{p}} C
  R(\tfrac{\la}t,A)^{k}x \right] \,\tfrac{d\la}{\la} \\
& & \quad + a \left[\left(\tfrac{1}t\right)^{k-\al-\sfrac{1}{p}}
   C(\tfrac{1}t+A)^{-k}\right].
\end{eqnarray*}
By our assumption \eqref{eq:W_C} and scaling invariance of
$\Ga$ and the measure $d\la/\la$ we obtain that $K(t)$ is uniformly
bounded, provided that the integral 
\[
\int_\Ga \bigl|f(\la)\bigr| \bigl|\la\bigr|^{-k+\al+\sfrac{1}{p}} \,d|\la|
\]
is finite. By the estimates in Lemma~\ref{lem:technical} we know
that $f\in O(|z|^{k-1-\eps})$ at zero and $f \in O(|z|^{-n})$ for any
$n\in\NN$ at infinity. Therefore the above integral is finite if
\[
k-1-\eps -k+\al+\sfrac{1}p >-1, \quad\text{i.e. if}  \quad \eps<\al+\sfrac{1}p.
\]
This however, due to our assumption on $\al$, may always
be satisfied by an appropriate choice of $\eps\in(0,1)$, and the proof
is done.
\end{proof}

To analyse $L^p$--admissibility for control operators,
assume that $\norm{T_{-k}(t)B}_{U\to X} \le M t^{-\ga}$ for some
$\ga\in\RR$. For $t>0$ fixed we thus have  
\begin{eqnarray*}
 \biggnorm{ \int_0^t T(t{-}s)Bu(s)\,ds }_X
 &\le& \int_0^t \bignorm{ T(t{-}s)B }_{U\to X} \, s^{-\al} \, \norm{ s^\al u(s) }_U\,ds\\
 &\le& c\int_0^t (t{-}s)^{-\ga} s^{-\al} \norm{ s^\al u(s) }_U\,ds.
\end{eqnarray*}
Let $k_{\al,\ga}(t,s)= \eins_{(0,t)}(s)(t{-}s)^{-\ga} s^{-\al}$ for $s,t\in(0,\tau)$.
Thus the study of the kernel $k_{\al,\ga}$ which may or may not induce
a bounded integral  operator $K_{\al,\ga}:L^p(0,\tau)\to
L^\infty(0,\tau)$ gives a sufficient criterion for (in)finite--time
$L^p$--admissibility of type $\al$ of control operators. 
\smallskip
For $p=1$, $K_{\al,\ga}$ is bounded if and only if $k_{\al,\ga}$ is
uniformly bounded, which, for finite $\tau$, is equivalent to
$\ga\le0$ and $\al\le0$. In this case we have
\[
  \norm{ K_{\al,\ga} }_{L^\infty(0,\tau)\to L^\infty(0,\tau)}
= \norm{ k_{\al,\ga} }_\infty 
=\frac{|\al|^{|\al|}|\ga|^{|\ga|}}{|\al+\ga|^{|\al+\ga|}}\,\tau^{|\al+\ga|}.
\]
For $\tau=\infty$, $k_{\al,\ga}$ is uniformly bounded if and only if $\al=\ga=0$.
\smallskip
For $p>1$ we use Hölder's inequality and obtain
\begin{eqnarray*}
      \int_0^t k_{\al,\ga}(t,s)|f(s)|\,ds
&\le& \biggl(\int_0^t \bigl(k_{\al,\ga}(t,s)\bigr)^{p'} \, ds
      \biggr)^{\sfrac{1}{p'}}\,\norm{ f }_p\\
      (\text{\rm subst. } s=t\si)
 &=&  \biggl(t^{1-(\ga+\al)p'} \int_0^1(1-\si)^{-\ga p'}\si^{-\al
        p'}\,d\si\biggr)^{\sfrac{1}{p'}}\,\norm{ f }_p.
\end{eqnarray*}
Convergence of the integral in the last line is equivalent to
$\ga<\sfrac{1}{p'}$ and $\al<\sfrac{1}{p'}$. Taking the sup over
$t\in(0,\tau)$ we see that, for finite $\tau$,
\[
 \norm{ K_{\al,\ga} }_{L^p(0,\tau)\to L^\infty(0,\tau)}\le c_{\al,\ga,p}\tau^{\sfrac{1}{p'}-(\ga+\al)} 
\]
provided that $\ga<\sfrac{1}{p'}$, $\al<\sfrac{1}{p'}$ and $\ga+\al\le
\sfrac{1}{p'}$, whereas, for $\tau=\infty$,
\[
 \norm{ K_{\al,\ga} }_{L^p(0,\infty)\to L^\infty(0,\infty)}\le c_{\al,\ga,p} 
\]
provided that $\ga+\sfrac{1}{p}<1$, $\al+\sfrac{1}{p}<1$ and
$\ga+\al+\sfrac{1}{p}=1$. Since, in fact (cf. \cite{Joergens:Int-Op}),
\[
 \norm{ K_{\al,\ga} }_{L^p(0,\tau)\to L^\infty(0,\tau)}
= \sup_{t\in(0,\tau)} \norm{ k_{\al,\ga}(t,\cdot) }_{L^{p'}(0,\tau)},
\]
we have just proven the following result.
\begin{proposition}\label{prop:bedingungen-fuer-control-zul-typ-al}
Assume that for some $\ga$, $\norm{T(t) B}_{U \to X} \le t^{-\ga}$, $t>0$. 
Then $K_{\al,\ga}$ is bounded $L^p(0,\tau)\to L^\infty(0,\tau)$ if and
only one of the following conditions holds:
\begin{equation}\label{eq:bedingungen-fuer-control-zul-typ-al}
\begin{array}{rlllll}
 (i)  & p=1 & \tau<\infty & \al\le0 & \ga\le0 & \\
 (ii) & p=1 & \tau=\infty & \al=0   & \ga=0   & \\
 (iii)& p>1 & \tau<\infty & \al+\sfrac{1}{p}<1 & \ga+\sfrac{1}{p}<1 
            & \al+\ga+\sfrac{1}{p}\le1\\
 (iv) & p>1 & \tau=\infty & \al+\sfrac{1}{p}<1 & \ga+\sfrac{1}{p}<1 
            & \al+\ga+\sfrac{1}{p}=1.
\end{array}
\end{equation}
Observe that condition $(iv)$ implies $\al>0$ and $\ga>0$. 
\end{proposition}

\begin{proof}[Proof of {Theorem~\ref{thm:control-Lp-zul-typ-alpha}}]
\ref{item:thm-control-zul-equiv-notw}.
Consider for $t>0$ the function $u(s) := \eins_{(\sfrac{t}2,t)}(s) u_0$.
Then 
\[
  \biggnorm{ \int_0^t T(t{-}s)B u(s) \,ds }_X
= \biggnorm{ \int_{\sfrac{t}2}^t A T(t{-}s)B u(s) \,ds }_{\dot{X}_{-1}}
= \bignorm{ \,[T(t)-T(\sfrac{t}2)]B u_0 }_{\dot{X}_{-1}}.
\]
We have $\norm{ u }_{L_\al^p([0,t], U)}  
=  \tfrac{1}{1+\al  p}\bigl(\sfrac{t}2\bigr)^{\al+\sfrac{1}p}$.
Let $f(z) := \exp(z){-}\exp(\sfrac{z}2)$. Notice that the function 
$z \mapsto z^{-\al-\sfrac{1}p} f(z) \in H_0^\infty(\Sec{\si})$ for all
$\si < \pihalbe$, so Theorem~\ref{thm:char-McI-Y-spaces} applies. We obtain 
\[
\cR(B) \subseteq \bigl( (\dot X(A_{-1}))_{-1}, (\dot X(A_{-1}))_1
                 \bigr)_{\einhalb(\al+\sfrac{1}p+1), \infty} 
=   \bigl( \dot{X}_{-2}, X \bigr)_{\einhalb(\al+\sfrac{1}p+1), \infty} 
=   \bigl( \dot{X}_{-k}, X \bigr)_{1-\theta, \infty}
\]
with $\theta = \sfrac{1}k\,(\sfrac{1}{p'}{-}\al)$, $k\in \NN$. The claim now 
follows from Theorem~\ref{thm:beschr-bedingung-und-interpol}.
\medskip
\ref{item:thm-control-zul-equiv-hinr}. 
By Theorem~\ref{thm:beschr-bedingung-und-interpol}, boundedness of $W_B$
is equivalent to $B: U \to (\dot{X}_{-k},X)_{\theta,\infty}$ with 
$\theta = 1 + \sfrac{1}k\, (\al{-}\sfrac{1}{p'})$. 
By analyticity of the semigroup this implies  
$\norm{ T(t)B } \le c t^{-\ga}$ with $\ga = \sfrac{1}{p'}-\al$
(cf. the arguments in the proof of
Lemma~\ref{lem:Abbildungsverhalten-von-analytischen-HG}).
Hence, if additionally to the assumptions of the theorem, 
$\al>0$ in case $p>1$ (or $p=1$ in case $\al=0$, respectively), 
Proposition~\ref{prop:bedingungen-fuer-control-zul-typ-al} gives the claim.
\medskip
\ref{item:thm-control-duale-bed-notw}.
Let (\ref{eq:def-Lp-zul-contr-falsch}) hold. For
$\Re(\la)>0$ and $u\in U$ we have
\[
(\lambda+A_{-k})^{-k} Bu = \frac{1}{(k-1)!}\, \int_{0}^{\infty}
t^{k-1} e^{-\lambda t} T_{-k}(t)B u\, dt.
\]
Then, by assumption
\begin{eqnarray*}
& &    \bignorm{ \Re(\la)^{k+\alpha-\sfrac{1}{p'}}\, (\lambda+A)^{-k-1} Bu } \\
&\le & \tfrac{1}{(k-1)!} \, \biggnorm{ \int_{0}^{\infty} t^{k-1}
       \Re(\la)^{k+\alpha-\sfrac{1}{p'}}\, e^{-\lambda t} T_{-k}(t)B
       u\, dt}\\
& = & \tfrac{1}{(k-1)!} \, \biggnorm{ \int_{0}^{\infty} 
       T_{-k}(t)B \bigl[t^{k-1} \Re(\la)^{k+\alpha-\sfrac{1}{p'}}\,
       e^{-\lambda t} \otimes u \bigr]\, dt}\\
&\le& \tfrac{K}{(k-1)!} \bignorm{h_\la(t) \otimes u }_{L^p_\al(\RR_+,U)}
 \le \widetilde K \norm{u}_U.
\end{eqnarray*}
Here, the uniform boundedness of the functions $h_\la(t) :=
\Re(\la)^{k+\alpha-\sfrac{1}{p'}} \, t^{k-1} e^{-\lambda t}$,
$\la>0$ in $L^p_\al(\RR_+)$ is shown similar to the proof of
Theorem~\ref{thm:obs-Lp-zul-typ-alpha}.
\medskip
\ref{item:thm-control-duale-bed-hinr}.
Without loss of generality we assume $k\ge 2$. We proceed in
two steps.\\
{\it Step 1:} For to show the existence of the integral in
\eqref{eq:def-Lp-zul-contr-falsch}, we chose some $x'\in (X_{-k})'$,
$x'\not=0$. Notice that $(X_{-k})'$ may be identified with $(X')_{k}$,
that is the domain $\cD((A')^k)$ with graph norm. We consider
\[
\int_0^\infty \bigl| \dual{ T_{-k}(t) B u(t)} {x'} \bigr|\,dt
=
\int_0^\infty \bigl| \dual{t^{-k+1} A^{-k+1} (tA)^{k-1} T_{-k}(t)B u(t)} {x'}
\bigr|\,dt.
\]
As we did in the proof of Theorem~\ref{thm:obs-Lp-zul-typ-alpha}, we
decompose $F_k(z) = z^{k-1} e^{-z}$ as $F_k(z) = \varphi(z)
\psi(z)$ with $\varphi$, $\psi$ as in \eqref{eq:decomp2}. 
We obtain 
\begin{eqnarray*}
&   & \int_0^\infty \bigl| \dual{T_{-k}(t) B u(t)} {x'} \bigr|\,dt \\
& = & \int_0^\infty \bigl| \dual{t^{-k+1} A_{-k}^{-k+1} \varphi(tA_{-k})
      \psi(tA_{-k}) B u(t)} {x'} \bigr|\,dt.
\eqnaintertext{
Notice that by sectoriality of $A$, the operators
$A_{-k}(\mu{+}A_{-k})^{-1}$, $\mu>0$ are uniformly bounded. 
Moreover, $\lim_{\mu\to 0+} A_{-k}(\mu{+}A_{-k})^{-1} Bu = Bu$ 
in $X_{-k}$ since $A_{-k}$ has dense range in $X$. 
Therefore, applying Fatou's lemma and writing $B_\mu :=
A_{-k}(\mu{+}A_{-k})^{-1} B$, we have
}
&\le& \liminf_{\mu\to 0+} \int_0^\infty \bigl| \dual{t^{-k+1}
      A_{-k}^{-k+1} \varphi(tA_{-k}) \psi(tA_{-k})
      B_\mu u(t)} {x'} \bigr|\,dt
\eqnaintertext{
Notice that $\psi(tA_{-k}) B u(t) \in \cD(A_{-k}^{-k+1+\eps})$ whereas 
$\psi(tA_{-k}) B_\mu u(t) \in \cD(A_{-k}^{-k+\eps})$. This observation
allows to interchange the operators $\varphi(tA_{-k})$ and
$A_{-k}^{-k+1}$ as follows
}
& = & \liminf_{\mu\to 0+} \int_0^\infty \bigl| \dual{\varphi(tA_{-k})
      t^{-k+1} A_{-k}^{-k+1} \psi(tA_{-k}) B_\mu
      u(t)} {x'} \bigr|\,dt \\
& = & \liminf_{\mu\to 0+} \int_0^\infty \bigl|
      \dual{t^{1+\sfrac{1}{p'}-k} A_{-k}^{-k} \psi(tA_{-k}) B_\mu u(t)}
      {t^{-\sfrac{1}{p'}}\varphi(tA_{-k})'x'} \bigr|\,dt\\
&\le& \liminf_{\mu\to 0+} \; \biggnorm{ t\mapsto
      \underset{=:L_\mu(t)}{\underbrace{\bigl[
       t^{\sfrac{1}{p'}-\al-k+1} A_{-k}^{-k-1}
      \psi(tA_{-k}) B_\mu \bigr]}} t^\al u(t)}_{L^p(\RR_+, X)}
      \bignorm{ \varphi(tA')x' }_{L^{p'}(\RR_+,dt/t,X')}.
\end{eqnarray*}
Notice that by assumption on $A'$, the $L^{p'}$--norm has an estimate
against the norm of $x'$, whence the existence of the integral is
proven if we show the uniform boundedness of the operators $L_\mu(t)$ for
$t>0$ and $\mu>0$. This step is very similar to the proof of
uniform boundedness of the family $K(t)$, $t>0$ in the proof of
Theorem~\ref{thm:obs-Lp-zul-typ-alpha}: we apply Lemma~\ref{lem:technical} with
$m=k-1$ to the function $\psi(z)$ and obtain for fixed $t>0$ 
\[
A_{-k}^{-k+1} \psi(tA_{-k}) B_\mu u(t)
=
\tfrac{(k-1)!}{2\pi i} \int_\Ga f(\la) t^{k-1} R(\la,tA_{-k})^{k} B_\mu u(t)
\,d\la
+ a t^{k-1} (I+tA_{-k})^{-k} B_\mu u(t).
\]
Now for $u \in U$ write
\begin{eqnarray*}
L_\mu(t)u
& = & \tfrac{(k-1)!}{2\pi i} \int_\Ga f(\la) t^{\sfrac{1}{p'}-\al}
       R(\la,tA_{-k})^{k} B_\mu u\,d\la 
  +    a t^{\sfrac{1}{p'}-\al} (I+tA_{-k})^{-k}  B_\mu u(t)\\
& = & \tfrac{(k-1)!}{2\pi i} \int_\Ga f(\la) \la^{-k-\al+\sfrac{1}{p'}}
       \left[A_{-k}(\mu{+}A_{-k})^{-1}\right]
       \left[\left(\tfrac{\la}t\right)^{k+\al-\sfrac{1}{p'}+1}
         R(\tfrac{\la}t,A_{-k})^{k} B\right]  u\,\tfrac{d\la}\la \\
&   &  \quad +   a \left[A_{-k}(\mu{+}A_{-k})^{-1}\right]
       \left[\left(\tfrac{1}t\right)^{k+\al-\sfrac{1}{p'}}
       (\tfrac{1}t+A_{-k})^{-k} B\right]  u.
\end{eqnarray*}
Therefore, by the assumption \eqref{eq:W_B} and a similar calculation to
\eqref{eq:beschraenktheit-zieht-sich-auf-Gamma-hoch} the set
$\{L_\mu(t):\; t>0, \mu>0\}$ is bounded in $B(U,X)$ provided that the
integral 
\[
\int_\Ga \bigl|f(\la)\bigr| \, |\la|^{-k-\al+\sfrac{1}{p'}} \,d|\la|
\]
is finite. Since $f \in O(|z|^{k-1-\eps})$ in zero and $f\in
O(|z|^{-n})$ for any $n\in\NN$, this boils down to
\[
k-1-\eps-k-\al+\sfrac{1}{p'} > -1, \quad \text{ i.e. to } \quad \al<\sfrac{1}{p'}-\eps,
\]
which, due to our assumption on $\al$, always may be satisfied by some 
$\eps\in(0,1)$.
\medskip
{\it Step 2:} Now let $x' \in X'$. We show that 
$t\mapsto \dual{T(t) Bu(t)}{x'} \in L^1(\RR_+)$ with a norm
estimate against $K' \norm{ u }_{L^p_\al(R_+,U)} \norm{x'}_{X'}$. To this
end we first notice that by analyticity of the semigroup, $T_{-k}(t) B
u(t) \in X$ for $t>0$. Moreover, for $t$ positive, 
$T_{-k}(t) B u(t) = \lim_{\eps\to 0} T_{-k}(t+\eps) Bu(t)$ in
$X$. Therefore, Fatou's lemma yields
\[
    \int_0^\infty \bigl| \dual{T_{-k}(t)Bu(t)}{x'}\bigr| \,dt
\le \liminf_{\eps\to 0} \int_0^\infty \bigl| 
       \dual{T_{-k}(t) Bu(t) } { T_{-k}(\eps)'x' }\bigr| \,dt.
\]
Notice that $y_\eps' := T_{-k}(\eps)'x' \in \cD((A')^k)$ and by step one, 
\[
      \int_0^\infty \bigl| \dual{T(t)Bu(t)} {y_\eps'}\bigr| \,dt
  \le K \norm{ u }_{L^p_\al} \norm{ y_\eps' }_{X'}.
\]
Since $\norm{ T_{-k}(\eps)'x'} \le K_0 \norm{x'}$ the integral
$\int_0^\infty T(t)Bu(t) \,dt$ exists as a Pettis integral in
$X$. The above argumentation shows that we have a bounded linear mapping
$\Phi: L^p(\RR_+,U) \to X''$, $\Phi(u) = \int_0^\infty T_{-k}(t)
B u(t)\,dt$. If $u$ is a step function with compact support that does
not contain zero, the integral in question even exists as a Bochner
integral. In this case, it takes values in $X$ by analyticity of
the semigroup $T(\cdot)$. Since such step functions are dense in $L^p(\RR_+, U)$ (recall
$p<\infty)$ we obtain $\cR{\Phi} \subseteq X$ and thus $\Phi$ is
necessarily bounded from $L^p(\RR_+, U)$ to $X$. This finishes our
proof.
\end{proof}
\subsection*{Regularity and Wellposedness}

\begin{proof}[Proof of {Lemma~\ref{lem:Abbildungsverhalten-von-analytischen-HG}}]
It is not hard to see, that for bounded analytic semigroups and
$k\in\NN$, $\norm{t^k A^k T(t)} \le c_k <\infty$. 
Indeed, by the elementary functional calculus for sectorial operators
(cf. \cite{McIntosh:H-infty-calc}) and substituting $tz =\la$ one has
\begin{eqnarray*}
\bignorm{ (tA)^k T(t)} 
&=& \biggnorm{ \tfrac{1}{2\pi i} \int_\Ga (tz)^k e^{-tz} R(z, A)\,dz}\\
&\le& M \int_\Ga |tz|^k e^{-t\,\text{Re}(z)} \,\tfrac{|dz|}{|z|}
= M \int_\Ga |\la|^k e^{-\text{ Re }\la} \,\tfrac{|d\la|}{|\la|} =: c_k
<\infty.
\end{eqnarray*}
This shows $\norm{T(t)}_{X \to \dot{X}_k} \le c_k t^{-k}$, $t>0$.
On the other hand, clearly $\norm{T(t)}_{X \to X} \le c_0$,
$t\ge 0$. By real interpolation we obtain immediately 
$\norm{ T(t)}_{X \to Z} \le k_1 t^{-\theta k}$.
Similarly, considering $\dot{X}_{-k}$ and $X$ in place of $X$ and
$\dot{X}_k$, one obtains $\norm{ T(t)}_{W \to X} \le k_2 t^{-\si k}$.
Both estimates together give the claim by the semigroup property.
\end{proof}
\begin{proof}[Proof of {Proposition~\ref{prop:PruessSimonett}}]
Let $\phi(s) := \bigl|(1+s)^\al -1\bigr|$. Then 
\[
\norm{ (Tf)(t) }_Y \le M \int_0^t \tfrac{1}{t-s} \phi(\tfrac{t-s}s)
\norm{f(s)}_U\,ds,
\]
whence $T$ is pointwise bounded by a multiple of the scalar integral operator
\[
(\widetilde T u)(t) := \int_0^t \tfrac{1}{t-s} \phi(\tfrac{t-s}s) u(s) \,ds
\]
which has the kernel $k(t,s) = \eins_{[0,t]}(s)\tfrac{1}{t-s}
\phi(\tfrac{t-s}s)$. Notice that, substituting $s=t\si$,
\[
\bignorm{k(t,\cdot)}_{p'}^{p'}
= 
\int_0^t \biggl| \frac{(1+\tfrac{t-s}s)^\al-1}{t-s}\biggr|^{p'}\,ds
=
t^{1-p'} \int_0^1 \biggl| \frac{\si^{-\al}-1}{\si-1}\biggr|^{p'}\,d\si
\le t^{1-p'} \widetilde c
\]
where $\widetilde c=\widetilde c(p,\al) <\infty$ since 
letting $g(\si) := \si^{-\al}$, the limit for $\si \to 1$ equals
$g'(1)=-\al$. It follows by Hölders inequality that 
$|\widetilde Tu(t)| \le c M \norm{u}_p t^{-\sfrac{1}p}$ (notice
$(1-p')/p'=-\sfrac{1}p$). Therefore, 
\[
    \la \, \mu\bigl( \{ t>0: \; |\widetilde Tu(t)|\ge \la\}\bigr)^{\sfrac{1}p} 
\le \la \, \mu\bigl(\{ t>0: \; M \norm{u}_p t^{-\sfrac{1}p}\ge \la\}\bigr)^{\sfrac{1}p} 
 =  M \norm{u}_p,
\]
showing that $\widetilde T$ is of weak type $(p,p)$ for 
every $p\in(1,\infty)$. By Marcinkiewicz interpolation, 
$\widetilde T$ is bounded on $L^p(\RR_+)$ for $p\in (1,\infty)$, which
implies the assertion.
\end{proof}

\begin{proof}[Proof of {Theorem~\ref{thm:PruessSimonett}}]
Let $\cF$ denote the convolution operator acting $L^p(\RR_+,U) \to L^p(\RR_+,Y)$ and
$\cF^\al$ denote the same operator acting $L_\al^p(\RR_+,U) \to L_\al^p(\RR_+,Y)$.
Further let $\Phi_\al: L_\al^p(\RR_+, \cdot) \to L^p(\RR_+, \cdot)$ be
the canonical isometric isomorphism given by $(\Phi_\al f)(t) := t^\al f(t)$.
Then $T := \Phi_\al \cF^\al \Phi_\al^{-1} - \cF$ satisfies
\[
(Tu)(t) = \int_0^t \bigl[\bigl(\tfrac{t}s\bigr)^\al - 1\bigr] K(t-s)\, u(s)\,ds
\]
with $K(t-s) = CT(t-s)B$. By analyticity of the semigroup, 
$K\in C(\RR_+, B(U,Y))$ and by hypothesis $\norm{K(t)} \le M/t$. Therefore,
Proposition~\ref{prop:PruessSimonett} applies and yields boundedness of $T:
L^p(\RR_+,U)\to L^p(\RR_+,Y)$. Hence $\cF$ is bounded
if and only if $\cF^\al$ is.
\end{proof}

\begin{proof}[Proof of {Theorem~\ref{thm:nonlin}}]
We let $v:=T(\cdot)x_0$ and denote, for $\rho,\tau>0$ to be fixed later,
\[
 \Sigma_{\rho,\tau}=\{x\in C([0,\tau],X)\cap L^p_\al([0,\tau],Z):
                      x(0)=x_0, \norm{x-w}_\Sigma\le\rho \}
\]
where $\norm{x}_\Sigma:=\max\{\norm{x}_{C([0,\tau],X)},\norm{x}_{L^p_\al([0,\tau],Z)})$.
Observe that we dropped $\tau$ in notation of the norm, and that 
$\Sigma_{\rho,\tau}$ is complete for the metric induced by $\norm{x}_\Sigma$.
We let $\Gamma x:=v+T(\cdot)*B(F(x)-F(x_0))Cx$ for $x\in\Si_{\rho,\tau}$
and shall choose $\rho$ and $\tau$ such that $\Gamma$ is a contraction on
$\Sigma_{\rho,\tau}$. Then Banach's fixed point theorem ends the proofs.
If $f$ is a simple function with values in $U$ with compact support in
$(0,\tau)$, then by analyticity of the semigroup $T(\cdot)B\ast f$ is a
continuous, $X$--valued function. Since such functions are dense
in $L^p_\al(U)$ (recall $p<\infty$), $T(\cdot)B\ast$ maps $L^p_\al(U)$
boundedly into $C([0,\tau], X)$. We let $c_v(\tau):=\norm{v-x_0}_{C([0,\tau],X)}$ and 
$l_v(\tau):=\norm{v}_{L^p_\al([0,\tau],Z)}$. In the following we shall drop the time 
interval in the norms. The assumptions imply that 
$\norm{T(\cdot)B*u}_\Sigma\le K\norm{u}_{L^p_\al(U)}$ for some $K>0$. We write $L$ for 
the Lipschitz-constant of the function $F$.
Then we have, for $x\in \Sigma_{\rho,\tau}$,
\begin{eqnarray*}
 \norm{\Gamma x-v}_\Sigma
 &\le& K\norm{(F(x)-F(x_0))Cx}_{L^p_\al(U)}\\
 &\le& KL\norm{ C }\norm{x-x_0}_{L^\infty(X)}\norm{x}_{L^p_\al(Z)}\\
 &\le& KL\norm{ C }(\norm{x-v}_{C(X)}+\norm{v-x_0}_{C(X)})
              (\norm{x-v}_{L^p_\al(Z)}+\norm{v}_{L^p_\al(Z)})\\
 &\le& KL\norm{ C }(\rho+c_v(\tau))(\rho+l_v(\tau)).
\end{eqnarray*}
Similarly we obtain, for $x,\tilde{x}\in\Sigma_{\rho,\tau}$,
\begin{eqnarray*}
 \norm{\Gamma x-\Gamma \tilde{x}}_\Sigma
 &\le& K\norm{(F(x)-F(x_0))Cx-(F(\tilde{x})-F(x_0))C\tilde{x}}_{L^p_\al(U)}\\
 &\le& K(\norm{(F(x)-F(x_0))C(x-\tilde{x})}_{L^p_\al(U)}
        +\norm{(F(x)-F(\tilde{x})C\tilde{x}}_{L^p_\al(U)}\\
 &\le& KL\norm{ C }(\norm{x-x_0}_{L^\infty(X)}\norm{x-\tilde{x}}_{L^p_\al(Z)}
              + \norm{x-\tilde{x}}_{C(X)}\norm{\tilde{x}}_{L^p_\al(Z)})\\
 &\le& KL\norm{ C }(\rho+c_v(\tau)+\rho+l_v(\tau))\norm{x-\tilde{x}}_\Sigma.
\end{eqnarray*}
Now we choose $\rho>0$ such that $\eta:=4KL\norm{ C }\rho<1$ and then $\tau>0$ such that
$\max\{c_v(\tau),l_v(\tau)\}\le\rho$. Thus we obtain
\[
 \norm{\Gamma x-v}_\Sigma\le\eta\rho<\rho\quad\mbox{and}\quad
 \norm{\Gamma x-\Gamma \tilde{x}}_\Sigma\le\eta\norm{x-\tilde{x}}_\Sigma
\]
for $x,\tilde{x}\in\Sigma_{\rho,\tau}$, as desired.
\end{proof}

\subsection*{$L^q$--spaces and Besov spaces}
\begin{proof}[Proof of {Lemma~\ref{lem:kaesekaestchen}}]
The definition of Besov and Triebel-Lizorkin spaces
together with Minkowski's inequality yield, for any $s\in \RR$, 
\[
\begin{array}{ll}
B^s_{q,p}(\Om) \emb F^s_{q,p}(\Om) & \text{provided that }q \ge p \text{  and} \cr
F^s_{q,p}(\Om) \emb B^s_{q,p}(\Om) & \text{provided that }q \le p.
\end{array}
\]
This will be used in the sequel. First we show the 'if'--part, that is
we show $L^q(\Om) \emb B^0_{q,p}(\Om)$ in the case 
$p \ge \max(2,q)$, that is $(\sfrac{1}q,\sfrac{1}p) \in$ I where area
I is as depicted below. By the above embeddings of Besov and
Triebel-Lizorkin spaces, $L^q(\Om)=F^0_{q,2}(\Om) \emb F^0_{q,p}
\emb B^0_{q,p}$. 
\begin{minipage}{6cm}
  \setlength{\unitlength}{1mm}
  \begin{picture}(50,60)(-5,-5)
     \put(0,0)      {\vector(1,0){50}}   
     \put(0,0)      {\vector(0,1){50}}   
     \put(40,0)     {\line(0,1){40}}   
     \put(0,40)     {\line(1,0){40}}   
     \put(0,0)      {\line(1,1){40}}   
     \put(0,20)     {\line(1,0){40}}   
     \put(24,8){I}   \put(31,24){II}   
     \put(14,31){III}\put(4,14){IV}
     
     \put(0,-5){$0$}\put(20,-5){$\sfrac{1}{2}$}
     \put(40,-5){$1$}\put(45,2){$\sfrac{1}{q}$}
     \put(-5,0){$0$}\put(-5,20){$\sfrac{1}{2}$}
     \put(-5,40){$1$}\put(2,45){$\sfrac{1}{p}$}
  \end{picture}
\end{minipage}
\newlength{\meinelaenge}\setlength{\meinelaenge}{\textwidth}
\addtolength{\meinelaenge}{-6.3cm}
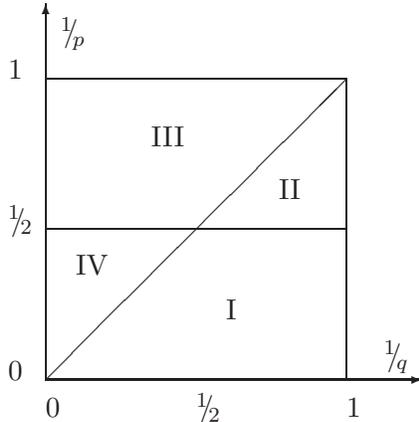
\begin{minipage}{\meinelaenge}\sloppy
Next, we consider the case $(\sfrac{1}q,\sfrac{1}p) \in$ III,
that is $p \le \min(q,2)$ and $(p,q)\not=(2,2)$. If $p \le 2 \le q$
and if $p<q$ we have
\begin{equation*}
B^0_{q,p}(\Om) \emb B^0_{q,2} \underset{\not=}{\emb} F^0_{q,2}(\Om) =
L^q(\Om),
\end{equation*}
and if $p,q<2$ and $p\le q$ we have
\begin{equation*}
B^0_{q,p}(\Om) \emb F^0_{q,p}(\Om)  \underset{\not=}{\emb}
F^0_{q,2}(\Om) = L^q(\Om).
\end{equation*}
Therefore, obviously $L^q(\Om) \not\emb  B^0_{q,p}(\Om)$.
\medskip
For counterexamples in area II and IV we construct 
specific functions $f\in L^q$ by wavelet decompositions
(cf. \cite{Meyer:ondelettes}), that show why the Besov norm 
cannot be estimated by the $L^q$--norm.
\end{minipage}
\medskip
Let $\Lambda$ be the set of all points
$\la=2^{-j}k+2^{-j-1}\eps$ 
where $j\in\ZZ$, $k\in\ZZ^n$ and $0\not=\eps\in \{0,1\}^n$.
Then every $\la\in\Lambda$ corresponds to unique $j$, $k$ and $\eps$.
Let $Q_\la$ be the dyadic cube defined by 
$Q_\la :=\{x\in \RR^n:\; 2^j x-k \in [0,1)^n \}$. Finally, by
\cite[Thm. III.8.1]{Meyer:ondelettes} chose some $1$-regular wavelet
basis $(\psi_\la)$ with compact support. Then  $\text{supp } \psi_\la
\subset c Q_\la$ for some $c>0$.  We let $\Lambda' := \{\la\in\Lambda:\;
\text{supp } \psi_\la \subset \Om\}$.

First we treat $(\sfrac{1}q,\sfrac{1}p)$ in area II, that is $q<p<2$.
By \cite[Thm VI.2.1]{Meyer:ondelettes}, for
$f=\sum_\la \al(\la) \psi_\la(x)$ in $L^q(\RR^n)$, we have
equivalence 
\begin{equation}
  \label{eq:eq-Lp-norm-in-waveletbasen}
\norm{f}_{L^q} 
\sim
\biggnorm{ \biggl(\sum_{\la\in \Lambda} |\al(\la)|^2 |Q_\la|^{-1}
  \eins_{Q_\la}\biggr)^\einhalb }_{L^q}.
\end{equation}
In the following, it will be sufficient to consider only functions $f$
that decompose in a finite sum. If $Q \subset \Om$ for some $Q =
Q_{\la_0}$, $\la_0\in\Lambda'$,  set $\al(\cdot)$ such that only
dyadic sub-cubes of $Q$ are  considered in the above summation: if
$Q_\la$ belongs to the $j^{\text{th}}$ dyadic subdivision of $Q$ then
let $\al(\la) := \al_j$, otherwise let $\al(\al) :=0$. Then the
expression on the right hand side of
\eqref{eq:eq-Lp-norm-in-waveletbasen} 
$L^q$--norm of $f$ simplifies to
\[
\biggnorm{\biggl(\sum_{j=0}^N \sum_{\Lambda_j'} \bigl|\al(\la)\bigr|^2
  \eins_{Q_\la} \biggr)^\einhalb}_{L^q}
= \bigl|Q\bigr|^{\sfrac{1}{q}} \biggl(\sum_{j=0}^N |\al_j|^2
  \biggr)^\einhalb.
\]
On the other hand side, an equivalent $B^0_{q,p}$--norm of 
$f = \sum_\la \al(\la) \psi_\la(x)$ is given by 
\begin{equation}
  \label{eq:eq-B0qp-norm-in-waveletbasen}
\norm{f}_{B^0_{q,p}} 
\sim
\biggl(\sum_{j=0}^N \biggl( \biggl(\sum_{\la\in \Lambda_j'}
\bigl|\al(\la)\bigr|^q\biggr)^{\sfrac{1}q}
2^{-nj(\sfrac{1}q-\sfrac{1}2)} \biggr)^p\biggr)^{\sfrac{1}p},
\end{equation}
see \cite[VI.10.5]{Meyer:ondelettes}. But, for $\la\in\Lambda_j'$ such
that $Q_\la\subset Q$, $|\al(\la)| = |\al_j| 2^{j\sfrac{n}q}$, whence 
\[
\norm{f}_{B^0_{q,p}} 
\sim 
\biggl(\sum_{j=0}^N |\al_j| 2^{nj \sfrac{p}2} \biggr)^{\sfrac{1}p}.
\]
Therefore, setting $\al_j := 2^{nj \sfrac{p}2}$ for $j=0,\ldots,N$
and letting $N\to\infty$ reveals that $L^q(\Om) \emb B^0_{q,p}$
implies  $p\ge 2$.
\medskip
Finally, consider the case IV, that is $q>2$ and $2<p<q$. 
If we set the wavelet coefficients $\al(\la)$ of $f$ in
\eqref{eq:eq-Lp-norm-in-waveletbasen} such that the cubes
in $\{ Q(\la): \;\al(\la)\not=0\}$ are piecewise disjoint, then 
\begin{eqnarray*}
\norm{f}_q
\sim
\biggnorm{ \biggl(\sum_{\la\in \Lambda} |\al(\la)|^2 |Q_\la|^{-1}
  \eins_{Q_\la}\biggr)^\einhalb }_{L^q}
&=& 
\biggnorm{ \sum_{\la\in \Lambda} |\al(\la)| |Q_\la|^{-\einhalb}
  \eins_{Q_\la} }_{L^q} \\
&=&
\biggl( \sum_{\la\in \Lambda} |\al(\la)|^q |Q_\la|^{-\sfrac{q}2}
  |Q_\la| \biggr)^\sfrac{1}{q} \\
&=&
\biggl( \sum_j \biggl(\biggl( \sum_{\la\in \Lambda_j}
|\al(\la)|^q\biggr)^{\sfrac{1}q} |Q_\la|^{\sfrac{1}q-\einhalb}
\biggr)^q \biggr)^\sfrac{1}{q}.
\end{eqnarray*}
On the other hand side, notice that $|Q_\la| = 2^{-nj}$.
Thus, comparing the $L^q$--norm of $f$ with the equivalent
$B^0_{q,p}$--norm given by \eqref{eq:eq-B0qp-norm-in-waveletbasen}, we
find $L^q(\Om) \emb B^0_{q,p}(\Om)$ requires $q\ge p$ contradicting
the assumption $p<q$.
\end{proof} 

\def\SUBMITTED{Submitted}
\def\TOAPPEAR{To appear in }
\def\PREPARATION{In preparation }

\end{document}